%% file: main.tex
\documentclass{article}

\pdfoutput=1

\usepackage[T1]{fontenc}
\usepackage{microtype}

\usepackage[margin=1.35in]{geometry}

\usepackage{amssymb}
\usepackage{amsthm}
\usepackage{amsmath}

\PassOptionsToPackage{hyphens}{url}\usepackage{hyperref}
    \hypersetup{
        pdfauthor={Steven Mattis, Kyle Robert Steffen, Troy Butler, Clint N. Dawson, Donald Estep},
        pdftitle={Learning Quantities of Interest from Dynamical Systems for Observation-Consistent Inversion},
        pdfsubject={Original research},
        pdfkeywords={Stochastic Inverse Problems; Uncertainty Quantification; Quantity of Interest; Observation-Consistent; Dynamical Systems},
        pdfproducer={Latex with hyperref},
        pdfcreator={pdflatex},
        colorlinks,
        linkcolor=blue,
        citecolor=blue,
        urlcolor=blue,
        final
    }
\usepackage{authblk}

\usepackage{amsfonts}
\usepackage{graphicx}
\graphicspath{{figures/}}
\usepackage{subfig}
\usepackage{xcolor}
\usepackage[lined,ruled,linesnumbered]{algorithm2e}
\usepackage{diagbox}

\newcounter{rmrk}[section]

\numberwithin{equation}{section}

\newtheorem*{ack*}{Acknowledgment}
\numberwithin{equation}{section}

\newcommand{\set}[1]{\left\{#1\right\}}

\newcommand{\ratio}{\mathcal{R}(Q(\lambda))}
\newcommand{\ratiocluster}{\mathcal{R}_k(Q_k(\lambda))}

\newcommand{\pspace}{\Lambda}
\newcommand{\dspace}{\mathcal{D}}

\newcommand{\cbinitabrv}{\text{init}} 

\newcommand{\cbupdateeriorabrv}{\text{update}} 

\newcommand{\initdens}{\pi^{{\cbinitabrv}}}
\newcommand{\updatedens}{\pi^{{\cbupdateeriorabrv}}}

\newcommand{\predictdens}{\pi^{Q(\text{init})}}

\newcommand{\predictdenscluster}{\pi^{Q_k(\text{init})}}

\newcommand{\obsdens}{\pi^\text{obs}}

\newcommand{\obsdenscluster}{\pi^\text{obs,k}}

\DeclareMathOperator*{\argmin}{arg\,min}

\title{Learning Quantities of Interest from Dynamical Systems for Observation-Consistent Inversion}

\author[1]{S. Mattis}
\author[2]{K.R. Steffen}
\author[3]{T. Butler}
\author[2,4]{C.N. Dawson}
\author[5,6]{D. Estep}

\affil[1]{Department of Statistics, Colorado State University}
\affil[2]{Oden Institute for Computational Engineering and Sciences, The University of Texas at Austin}
\affil[3]{Department of Mathematical and Statistical Sciences, University of Colorado Denver}
\affil[4]{Department of Aerospace Engineering and Engineering Mechanics, The University of Texas at Austin}
\affil[5]{Canadian Statistical Sciences Institute}
\affil[6]{Department of Statistics and Actuarial Science, Simon Fraser University}

\begin{document}

\maketitle

\begin{abstract}
Dynamical systems arise in a wide variety of mathematical models from the physical, engineering, life, and social sciences.
A common challenge is to quantify uncertainties on model inputs (i.e., parameters) that correspond to a quantitative characterization of uncertainties on observable Quantities of Interest (QoI).
To this end, we consider a stochastic inverse problem (SIP) with a solution described by a pullback probability measure.
This is referred to as an observation-consistent solution since its subsequent push-forward through the QoI map matches the observed probability distribution on model outputs. 
A distinction is made between QoI useful for solving the SIP and arbitrary model output data. 
In dynamical systems, model output data are often given as a series of state variable responses recorded over a particular time window.
Consequently, the dimension of output data can easily exceed $\mathcal{O}(1E4)$ or more due to the frequency of observations, and the correct choice or construction of a QoI from this data is not self-evident.
We present a new framework, Learning Uncertain Quantities (LUQ), that facilitates the tractable solution of SIPs for dynamical systems.
Given ensembles of predicted (simulated) time series and (noisy) observed data, LUQ provides routines for filtering data, learning the underlying dynamics in an unsupervised manner, classifying the observations, and performing feature extraction to learn the QoI map.
Subsequently, time series data are transformed into samples coming from the underlying predicted and observed distributions associated with the QoI so that solutions to the SIP are computable. 
Following the introduction and demonstration of LUQ, numerical results from several SIPs are presented for a variety of dynamical systems arising in the life and physical sciences.
In the interest of scientific reproducibility, we provide links to our Python implementation of LUQ, as well as all data and scripts required to reproduce the results in this manuscript.

\end{abstract}


\input{Introduction}

\input{LUQ-framework}

\input{LUQ-implementation}

\input{observation-consistent}

\input{Numerics}

\input{Conclusions}

\section{Acknowledgments}

Troy Butler is supported by the National Science Foundation (DMS-1818847).
Clint Dawson's and Kyle Robert Steffen's work is supported in part by the National Science Foundation (DMS-1818941).
Donald Estep's work is undertaken, in part, thanks to funding from the Canada Research Chairs Program, and is partially supported by the National Science Foundation under grants DMS-1821210, DMS-1818777, and DMS-1720473, by Riverside Research under contract RADIAEM.IDIQ.05 PO\#00133, and by grants from the Natural Sciences and Engineering Research Council of Canada.
Steven Mattis's work is supported by the National Science Foundation (DMS-1818777).
The authors acknowledge the Texas Advanced Computing Center (TACC) at The University of Texas at Austin for providing HPC resources that have contributed to the research results reported within this paper.

 \appendix

\input{Software-dependencies}




\bibliographystyle{elsarticle-num}
\bibliography{references.bib}







\end{document}

%% file: Introduction.tex
\section{Introduction}\label{sec:intro}


In order for stakeholders to make knowledgeable, data-informed decisions that incorporate computational model simulations of dynamical systems, it is essential to provide useful quantitative characterizations of uncertainties on the spaces defined by model inputs and outputs. 
Broadly speaking, a forward uncertainty quantification (UQ) analysis studies how uncertainty on model inputs propagates to model outputs through a Quantity of Interest (QoI) map while an inverse UQ analysis studies how uncertainty on model outputs corresponds to uncertainty on model inputs.
The last several decades has seen the uncertainty quantification (UQ) community develop many formulations and solution methodologies for various forms of forward and inverse UQ analyses that are of interest to the broader scientific community, e.g., see \cite{Berger_TEST_1994, Wikle1998, Kennedy_O_JRSSSB_2001, Stuart_IP_2010, Sargsyan_NG_IJCK_2015, HASSELMAN20082596, WU201559, ROY20112131} and the references therein. 

The stochastic inverse problem (SIP) considered in this work involves the computation of pullback measures\footnote{Let $Q:\pspace\to\dspace$ denotes the QoI map between the space of model inputs, denoted by $\pspace$, and the space of model outputs, denoted by $\dspace$. If $P_\pspace$ is a probability measure on $\pspace$, then its push-forward on $\dspace$, denoted by $P_\dspace^Q$, is defined by $P_\dspace^Q(A)=P_\pspace(Q^{-1}(A))$ for all events $A\subset\dspace$. Here, $Q^{-1}(A)$ denotes the pre-image of $A$ under $Q$. On the other hand, if $P_\dspace$ is a probability measure on $\dspace$, then a pullback on $\pspace$ is any measure $P_\pspace$ satisfying the property $P_\dspace^Q(A)=P_\dspace(A)$ for all events $A\subset\dspace$.} on model inputs associated with a probability measure observed on model outputs. 
We refer to these pullback measures as observation-consistent solutions since they induce a push-forward measure through the QoI map that yields the observed probability measure. 
The formulation of this SIP and observation-consistent solutions are rooted in rigorous measure theory that directly handles the set-valued inverses common in QoI maps (e.g., see \cite{BBE11, BES12, BET+14, BJW18a}).
In this approach, the QoI map is not regularized, which separates this SIP and its solution from other inverse problem formulations and their solutions.
For example, in what is commonly referred to in the UQ community as the ``Bayesian inverse problem,'' the map is implicitly regularized through the choice of a prior in a Bayesian formulation, e.g., see \cite{0266-5611-24-5-055012, 0266-5611-26-2-025002, doi:10.1137/130933381}.
Previous analysis of the SIP and observation-consistent solutions involved an a priori specification of the QoI map \cite{BJW18a, BBE11, BES12, BET+14}.
In \cite{WWJ17}, an experimental design paradigm for the SIP chooses an optimal QoI map from a specified list of QoI for which data may be collected.
More recently, in \cite{BH20}, QoI maps are constructed using a principal component analysis of the eigenvectors of the Laplacian operator.
However, no previous work has studied the construction of a QoI map directly from temporal data in order to solve the SIP. 
A major contribution of this work is the development of the Learning Uncertain Quantities (LUQ) framework for constructing QoI maps from noisy data for a dynamical system. 
This framework is encoded within an open-source Python package \verb|LUQ| \cite{mattis_luq}.
We demonstrate this framework and the \verb|LUQ| software package on several dynamical system models of various conceptual and computational complexity, including models with Hopf bifurcations or shocks (discontinuous solutions) and a quasi-operational model of storm surge (coastal flooding) from hurricanes and cyclones. 

We first describe the interplay between the SIP and LUQ framework in the context of a classic mass-spring system 
with parameters related to physical properties modeled by spring constants (i.e., the ``rigidity'' of the system) and the amount of energy dissipated by the system (i.e., the internal friction of the system). 
Uncertainty arising from imperfections in manufacturing and design processes are modeled using probability distributions.
The effect of these uncertainties are observable indirectly through the analysis of system responses to various external stimuli.
The LUQ framework transforms the data of system responses into QoI samples enabling the construction of an observed probability distribution on the QoI. 
Solving the SIP then produces probability distributions on parameters that are consistent with this observed distribution on the QoI. 

The problem of constructing an observation-consistent distribution is important to the broad scientific and engineering community.
For instance, in the biological sciences, nonlinear systems of ordinary differential equations are often used to model many types of phenomena ranging from competition of species (e.g., the Lotka--Volterra equations), spread of disease (e.g., SIR models), and metabolic processes of living organisms (e.g., the Sel'kov model for glycolysis). 
The coefficients in these models are generally uncertain and modeled using probability distributions. 
Likewise, in coastal engineering applications, the modeling of waves or storm surge depends upon a variety of uncertain model inputs, such as initial conditions, meteorological forcing, and as coefficients determining the free surface and bottom stress parameterizations. 

%
%

All of these problems lead to the same SIP: Determine a distribution on model inputs that yields the observed distribution on the QoI.
However, in determining the QoI for dynamical systems using temporal observations, one faces several challenges.
First, while relatively little data may be collected spatially at each observation time, individual measurement devices can easily produce $\mathcal{O}(10^4)$ (or more) time steps for which data are collected, though the dynamical behavior may be characterized in a ``low-dimensional'' way. 
The problem is to perform dimension reduction by transforming the data sets into low-dimensional QoI.
However, this leads to a second technical challenge: A dynamical system may present strikingly different qualitative behaviors as model inputs are varied (e.g., due to bifurcations) or over different windows of time (e.g., due to transient or equilibrium behavior). 
The LUQ framework addresses both of these problems. 
We utilize unsupervised learning approaches applied to data that classify the different types of dynamics that occur.
This is followed by the training and optimal selection of classifiers on this labeled data so that observed temporal data are appropriately ``binned'' by their prevalent dynamics. 
To each type of learned dynamical system response, we subsequently construct different forms of the QoI maps that best describe the low-dimensional nature of the dynamics.


The rest of this manuscript is organized as follows. 
In Section~\ref{sec:LUQ}, we introduce the new conceptual framework Learning Uncertain Quantities (LUQ).
This is followed in Section~\ref{sec:LUQ-software} by an in-depth overview of the computational package \verb|LUQ| that is used to generate the results presented in this paper.
Section~\ref{sec:inversion} summarizes the application of outputs of LUQ to construct observation-consistent solutions. 
Section~\ref{sec:numerics} includes numerical results with applications from the life and physical sciences.
Finally, in Section~\ref{sec:conclusions}, concluding remarks along with on-going and future research directions are given.
In the interest of scientific reproducibility, \ref{app:dependencies} contains links to the \verb|LUQ| repository as well as the data sets and Python scripts used to generate the various figures and table data appearing in this manuscript.

%% file: LUQ-framework.tex
\section{Learning Uncertain Quantities (LUQ): Conceptual Framework}\label{sec:LUQ}



To help motivate and illustrate the various steps in this framework, we consider the equation for the amplitude of displacement, $y$, in a (damped) harmonic oscillator 
given by
\[
	y''(t) + 2cy'(t) + \omega_0^2 y(t) = f(t)
\]
where $c$ is interpreted as a damping constant and $\omega_0$ is interpreted as the natural frequency.
For simplicity, we assume that there is no external forcing (i.e., $f(t)=0$).
Any observed motion is determined by an initial displacement given by $y(0)=3$ and $y^\prime(0)=0$.
The dynamics of this system are qualitatively described as under-damped, critically damped, or over-damped, depending on the relationship between $c$ and $\omega_0$.
We ignore units except to interpret a single unit of time as a second.

To motivate the associated inverse problem, suppose this model is used to describe the motion of a manufactured mass-spring-dashpot system.
Due to imperfections in the manufacturing process, the parameters $c$ and $\omega_0$ possess aleatoric uncertainty (i.e., irreducible variability in their values) across the family of devices.
Purely for illustrative purposes, we assume that physically plausible bounds for the parameters are
\[
	0.1\leq c \leq 1 \quad \text{ and } \quad 0.5\leq \omega_0 \leq 1,
\]
and that the data-generating distributions for these parameters are independent Beta$(2,2)$ distributions over their respective domains.

Assuming that actual observations are polluted by measurement error given by an additive noise model, the observed data, denoted by $y_\text{obs}(t)$, are written as
\[
	y_\text{obs}(t) = y(t) + \eta(t)
\]
where at each observed time, the measurement error, denoted by $\eta(t)$, is assumed independent and identically distributed (i.i.d.) according to an $N(0,\sigma^2)$ distribution.
Purely for illustrative purposes, we set $\sigma=0.25$ and assume that 
data are collected at a rate of 100 Hz starting at $1$ second and ending at $6$ seconds for a total of 501 measurements.
In Section~\ref{sec:harmonic-times}, we discuss the impact of taking observations in a different window of time.

We simulate $300$ i.i.d.~experiments involving parameter values drawn from the Beta distributions for each parameter to generate an associated family of $300$ (noisy) observed data vectors, each of length $501$.
Now, we suppose that the data-generating distribution is unknown, and we attempt to compute a distribution on the parameters that is observation-consistent.
We choose  an ``initial'' distribution of independent uniform distributions over each parameter domain.
We simulate observation data using $2000$ i.i.d.~samples drawn from this initial distribution.
This data is referred to as ``predicted'' data that comes from the push-forward measure for the initial distribution.
As discussed in Section~\ref{sec:inversion}, we can use the subsequent samples of QoI data obtained from noisy and predicted data to construct density estimates and perform an analog of ``accept-reject'' to estimate an ``update'' to the initial distribution that is observation-consistent. 

\subsection{Filtering data (approximating dynamics)}\label{ssec:Cleaning}
\begin{figure}[htbp]
	\centering
	\includegraphics[width=0.4\textwidth]{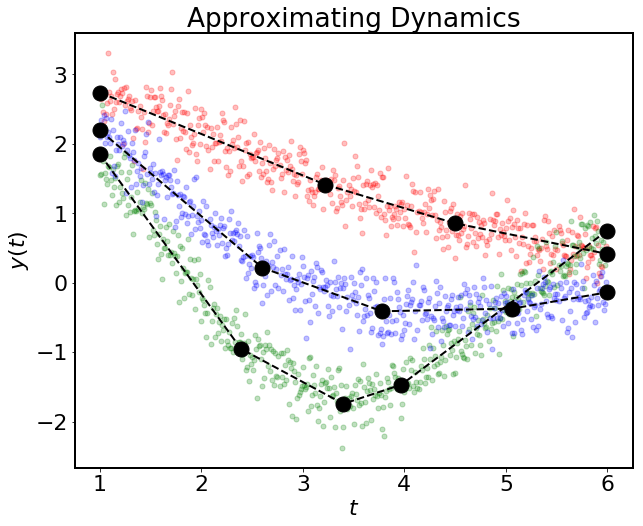}
	\includegraphics[width=0.4\textwidth]{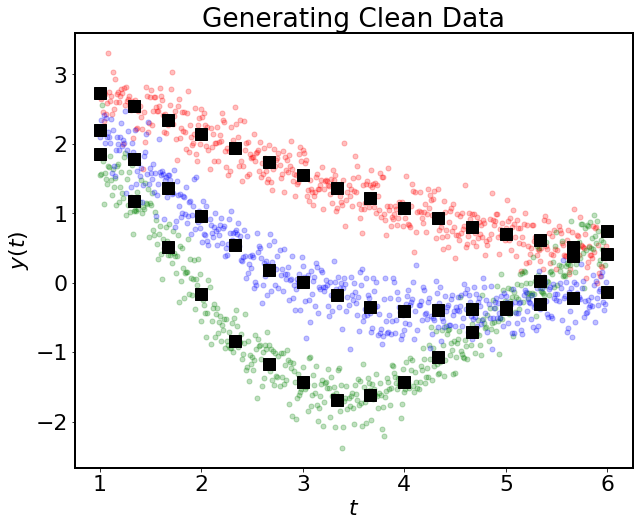}
	\caption{The small colored dots in each plot represent three different sets of noisy data filtered by a two-step process.
	Step 1 (left): approximate dynamics from noisy data using splines with number of knots and knot locations optimized according to Algorithm~\ref{alg:cleaning} in Section~\ref{sec:LUQ-software}.
	Step 2 (right): sample splines to generate filtered data.
	In the left plot, the dashed lines represent the splines constructed from each noisy time series, and the black disks denote the knots of each spline. 
	Observe that splines with the same number of knots may still have different knot locations due to the optimization algorithm.
	In the right plot, the regularly spaced black squares in each cloud of noisy data represent the filtered data sampled from the splines shown in the left plot.}
	\label{fig:harmonic-oscillator-cleaning-data}
\end{figure}

In the problems considered in this work, data are observed (or predicted) from a dynamical response of a physics-based model at some frequency over finite time intervals.
We produce filtered data by sampling, at potentially different frequencies,  approximations to the underlying dynamical response associated with each time series of data.
Below, we describe the basic procedure and give more specific implementation details in Section~\ref{sec:LUQ-software}.

Inspired by the work of \cite{UFW17} for adaptively constructing splines to solve linear inverse problems associated with potentially discontinuous functions, we use piecewise linear splines with both adaptive numbers of knots and adaptive knot placement to approximate underlying dynamical responses.
As described in \cite{UFW17}, it is then possible to approximate the underlying dynamical response (assuming a particular regularity or finite number of discontinuities) to arbitrary pointwise accuracy if both a sufficiently high frequency for collecting data and number of knots are used.

In Figure~\ref{fig:harmonic-oscillator-cleaning-data}, we illustrate this process of filtering noisy observations from three distinct dynamical system responses for the harmonic oscillator.
The left plot demonstrates the first step of approximating the underlying dynamics with adaptively-generated splines.
We observe that the different sets of noisy data require either different numbers of knots or knot placement in order to capture the underlying variability of the signal responsible for the noisy data.
To demonstrate the different uses and requirements of the ``noisy'' and ``filtered'' data in the analysis, the filtered data are generated at a lower frequency from the splines.
This is shown in the right plot.
Below, we ultimately observe that not as many filtered data are required to generate useful QoI whereas many noisy data may be required to generate accurate approximations of the dynamics using splines.

To ensure compatibility between the predicted and observed data sets, we also use predicted (noise-free) data to construct and sample splines to form ``filtered'' predicted data.
This may, in fact, be necessary in some scenarios such as when predictions come from controlled experiments containing measurement noise or if the measurement noise model is utilized when forming predictions.


\subsection{Clustering and classifying data (learning and classifying dynamics)}

%
%
The goal is to construct low-dimensional QoI that characterize dynamical behaviors to construct observation-consistent solutions on parameters. 
Since different dynamical behaviors may be characterized by different QoI, the first goal is to use (filtered) predicted data to determine the equivalent classes of dynamical behavior present in the (filtered) observed data.
This requires labeling the dynamics present in the predicted data set.
Clustering algorithms are a type of unsupervised learning algorithm that label data vectors using a metric to gauge the distance of a vector from the proposed ``center'' of the cluster (see \cite{xu2015comprehensive} for a comprehensive review of clustering algorithms).

To illustrate a particular clustering algorithm, we apply a $k$-means algorithm \cite{arthur2006k} on the predicted data for the harmonic oscillator.
The $k$-means algorithm is a centroid-based clustering algorithm.
Since the \verb|LUQ| module utilizes the \verb|scikit-learn| Python library, other clustering algorithms are readily available within LUQ as discussed in Section~\ref{sec:LUQ-software} below.
In this case, we have a priori knowledge of three potential types of dynamics in the harmonic oscillator, so we propose three centroids in the $k$-means algorithm.
A future work will consider approaches for adaptively choosing an optimal number of clusters when there is no a priori knowledge of dynamical behaviors present in an application.
\begin{figure}[htbp]
	\centering
	\includegraphics[width=0.8\textwidth]{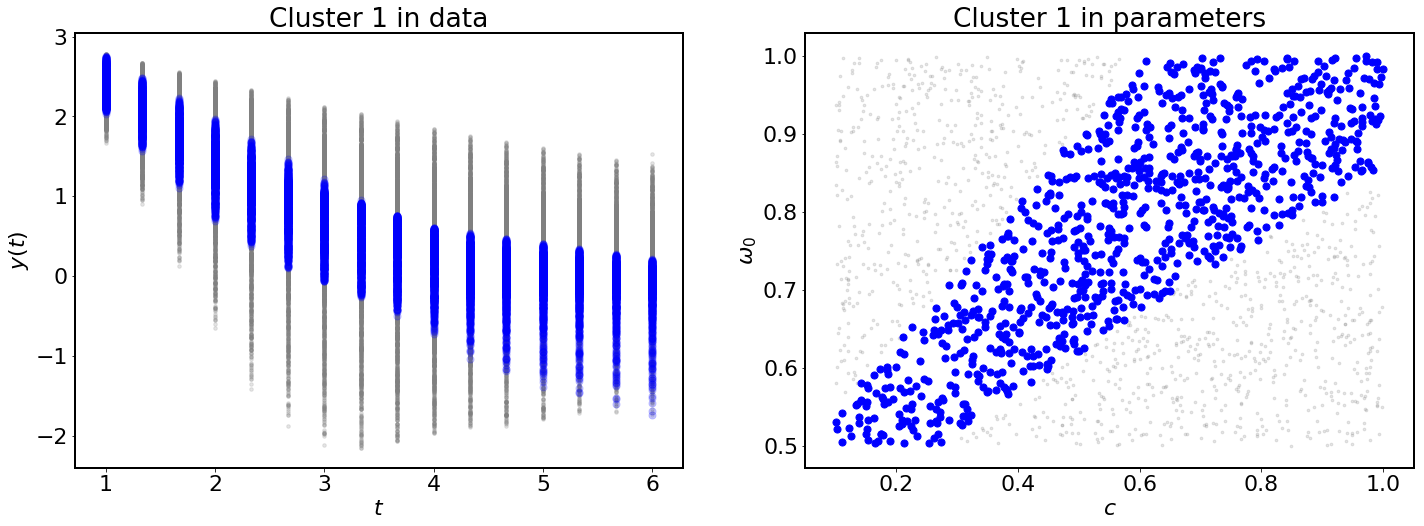}\\
	\includegraphics[width=0.8\textwidth]{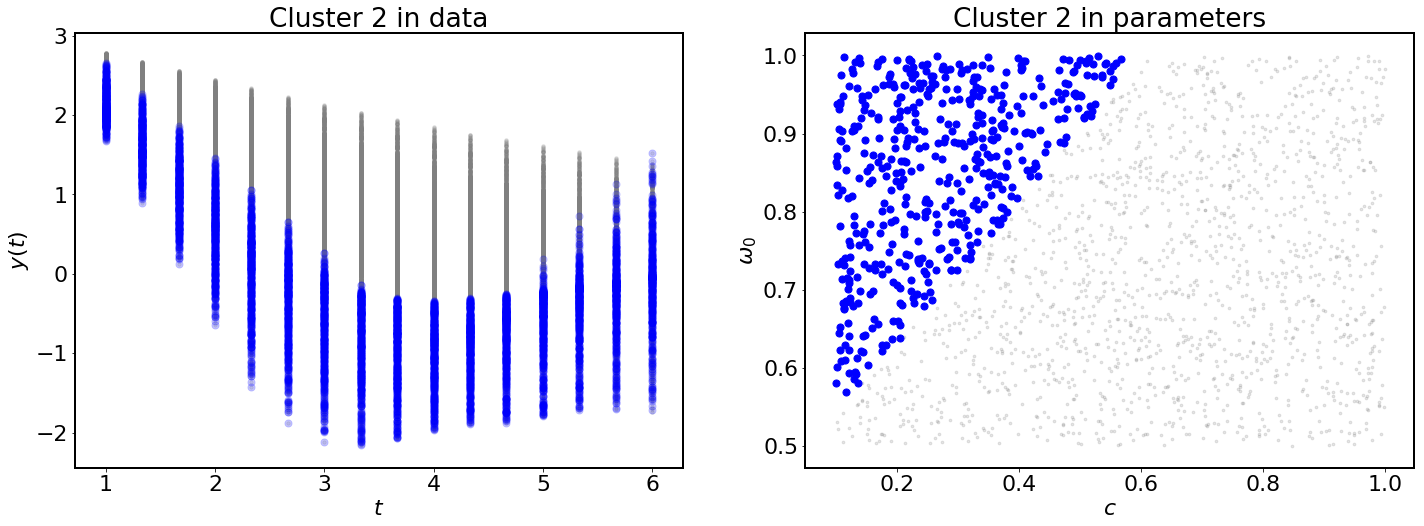}\\
	\includegraphics[width=0.8\textwidth]{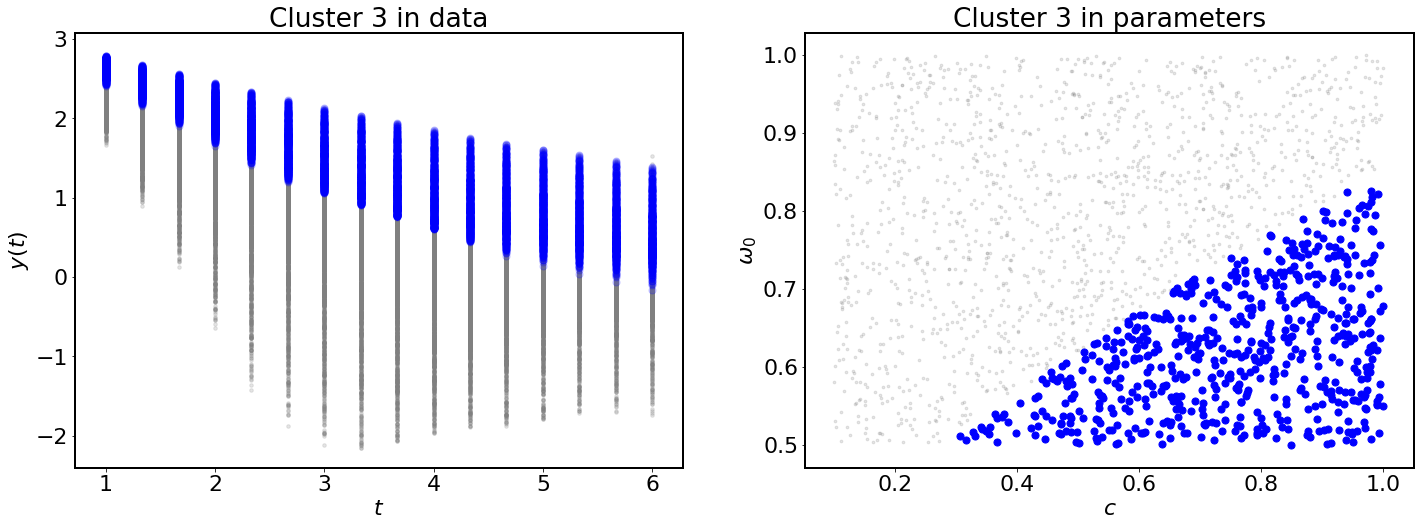}
	\caption{The large blue dots illustrate the cluster in each plot (in time on the left, in the parameter space on the right) using unsupervised learning in the form of $k$-means clustering. The smaller gray dots belong to different clusters. From top to bottom, we show clusters 1, 2, and 3, respectively.}
	\label{fig:harmonic-oscillator-clustering}
\end{figure}
Figure~\ref{fig:harmonic-oscillator-clustering} illustrates the output of this clustering.
The left column of plots show the clustering has grouped similar dynamical responses in time.
The right column of plots show the clustering of samples in parameter space inferred from the clustering of filtered data.



\begin{figure}[htbp]
	\centering
	\includegraphics[width=0.32\textwidth]{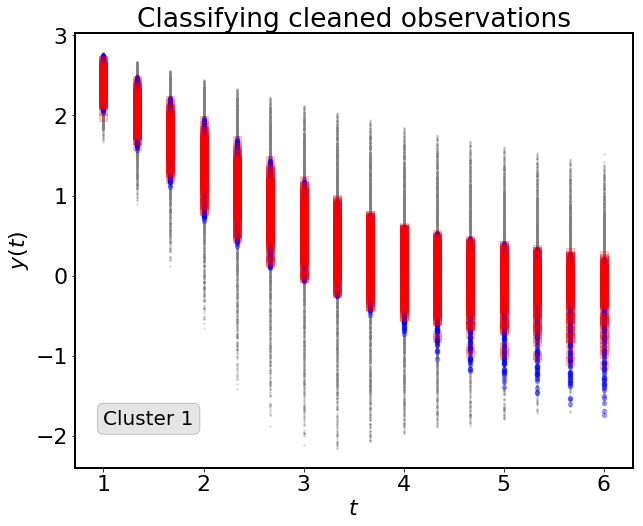}
	\includegraphics[width=0.32\textwidth]{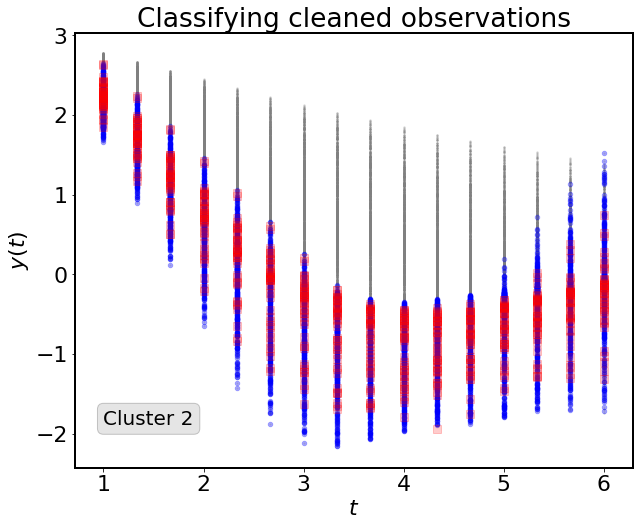}
	\includegraphics[width=0.32\textwidth]{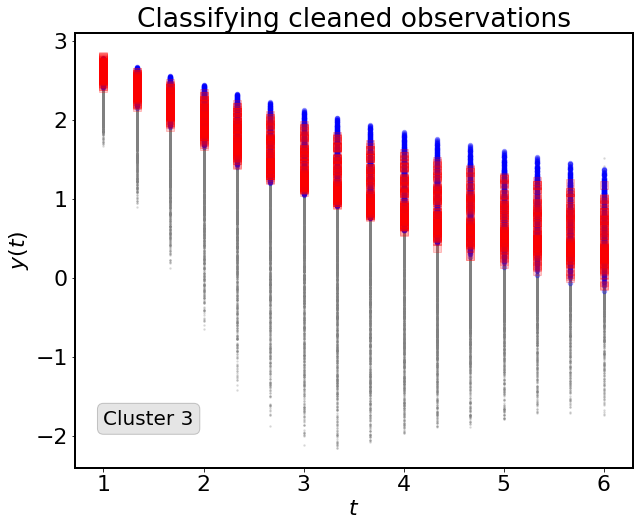}
	\caption{The large red dots illustrate the classified (filtered) observed data associated with a cluster (blue dots) learned from the set of all (filtered) predicted data (small gray dots). From left-to-right illustrates the classification of observed data into clusters 1, 2, and 3, respectively.}
	\label{fig:harmonic-oscillator-classifying-data}
\end{figure}
A classifier is a type of supervised learning algorithm (see \cite{kotsiantis2007supervised,alpaydin2014introduction}) that uses labeled training data to tune the various classifier parameters (not to be confused with the model parameters) so that non-labeled observed data can be properly labeled (i.e., classified).
At a high level, classifiers usually partition the training data into two subsets of data: one used to tune the classifier parameters and the other to test the quality of these tuned parameters by tracking the rate of misclassification.
This is referred to as cross-validation and is also useful in avoiding over-fitting the classifier (i.e., over-tuning the classifier parameters) to the entire set of training data \cite{kim2009estimating,hastie2009elements}.

A typical workflow for optimizing a classifier is to randomly split the training data into two subsets, perform cross-validation, and repeat this process some predetermined number of times to choose the classifier parameters that give the lowest misclassification rate \cite{hastie2009elements}.
The number of times, $k$, that the process is repeated is referred to as performing $k$-fold cross-validation.
It is also sometimes necessary to apply the so-called ``kernel trick'' to the data to transform it (often nonlinearly) into a higher-dimensional space where the construction of a classifier is made simpler \cite{boser1992training,scholkopf2001kernel}.
In this work, we use the labeled predicted data from the clustering step as the training data and restrict focus to kernel-based support vector machines (SVMs) \cite{cortes1995support}.
This is summarized in more detail in Section~\ref{sec:LUQ-software}.

For the harmonic oscillator, a linear kernel SVM is trained on the labeled predicted data using ten-fold cross-validation, leading to a classifier with a misclassification rate of approximately 0.25\%. 
Figure~\ref{fig:harmonic-oscillator-classifying-data} shows the classification of observed data using this kernel SVM.
The large red dots are the classified observed time series that are plotted on top of the associated predicted cluster data shown as blue dots.
The gray dots show the range of all predicted data.
That the red dots appear to be contained within the vertical ranges of the blue dots at all times indicates the SVM is properly classifying the dynamics present in the observed data sets. 

\subsection{Feature extraction (learning quantities of interest)}

\begin{figure}[htbp]
	\centering
	\includegraphics[width=0.32\textwidth]{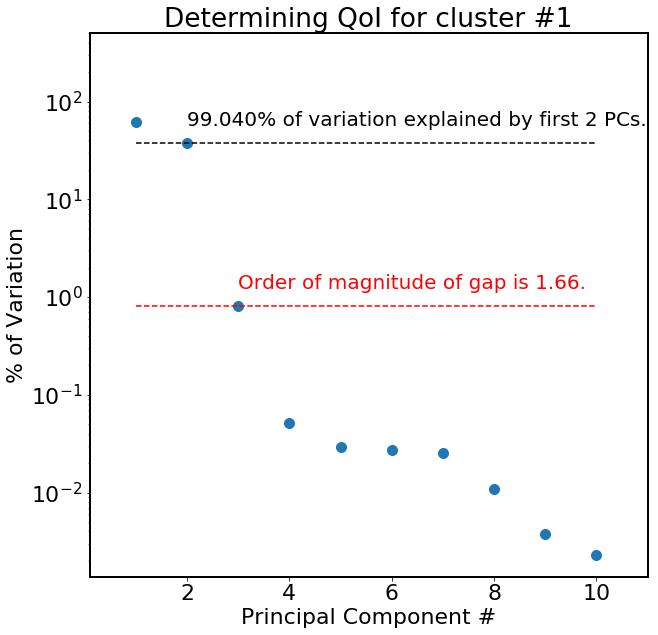}
	\includegraphics[width=0.32\textwidth]{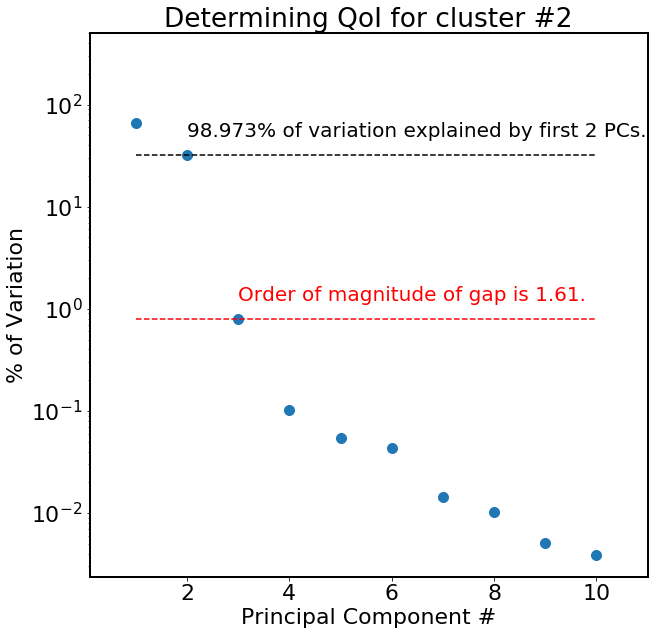}
	\includegraphics[width=0.32\textwidth]{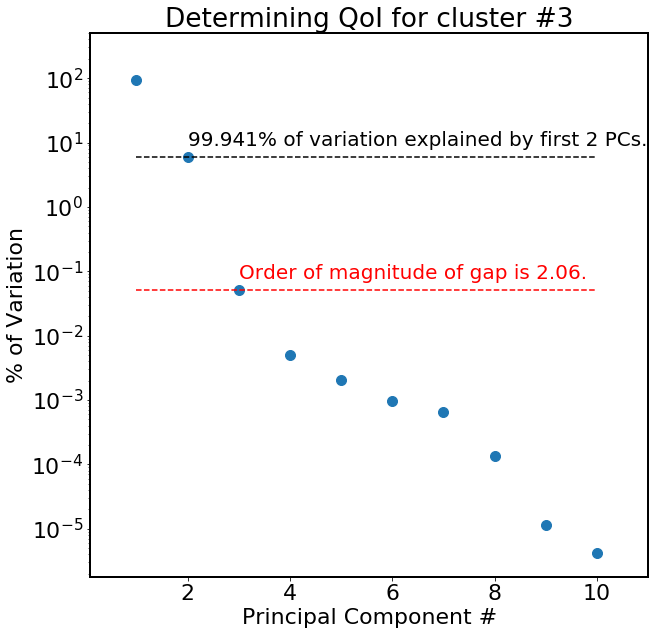}
	\caption{Performing kernel-based PCA on each cluster of predicted data to learn the QoI. From left to right is the output for clusters 1, 2, and 3, respectively. The amount of variation in the data summarized by the first two principal components is given above the black (upper) dashed line. The magnitude of the spectral gap to the third principal component is given above the red (lower) dashed line.}
	\label{fig:harmonic-oscillator-learning-QoI}
\end{figure}
Feature extraction algorithms generally attempt to reduce the dimension of a data space into a relatively small number of quantities that explain most of the variation observed in the data \cite{guyon2008feature,alpaydin2014introduction}.
In this work, we use kernel-based principal component analysis (PCA) \cite{pearson1901liii,scholkopf1997kernel,mika1999kernel,abdi2010principal}, which is one of the most popular approaches to performing feature extraction.
Kernels are measures of similarity, i.e., $s(a, b) > s(a, c)$ if objects $a$ and $b$ are considered ``more similar'' than objects $a$ and $c$ \cite{zhang2007local}.
At a high level, we use a kernel to transform the predicted data of each cluster into a space where a standard PCA is applied.
The percentage of variation explained by the first few dominant principal components is then computed.
This is carried out separately on each cluster, and then the dominant principal components are subsequently used as the QoI associated with the dynamics present in a particular cluster.

We seek to determine a kernel-based PCA such that nearly all of the variation in the data in a particular cluster is explained by the first $p$ principal components where $p$ is the dimension of the parameter space.
This is due to the basic fact that the dynamics in each cluster are parameterized by the $p$ parameters, and the goal is to update the distribution for {\em all} of the $p$ parameters.
The QoI can be thought of as describing a change of basis from the $p$ parameters into the space of dynamics described by each cluster.
In Section~\ref{sec:LUQ-software}, more details are provided about the various kernels applied to the data within \verb|LUQ| and the processes available for choosing the best performing kernel for the PCA.

In Figure~\ref{fig:harmonic-oscillator-learning-QoI}, we visualize the results of a linear kernel PCA in terms of the percentage of variation explained by the first two principal components as well as the magnitude of the ``spectral gap'' between the second and third principal components.
Here, we truncate the plots of the principal components at the maximum number of knots allowed in the adaptive splines used to filter the data.
Note that more than 98\% of the variation in the predicted time series data in each cluster is explained by just two principal components.

%% file: LUQ-implementation.tex
\section{Learning Uncertain Quantities: Implementation}\label{sec:LUQ-software}
The \verb|LUQ| Python package \cite{mattis_luq} provides simple implementations of the algorithms for learning uncertain quantities outlined in Section~\ref{sec:LUQ}.
\verb|LUQ| utilizes several publicly available Python packages that are commonly used for scientific computing (\verb|NumPy| \cite{oliphant2006guide} and \verb|SciPy| \cite{2020SciPy-NMeth}) and machine learning (\verb|scikit-learn| \cite{scikit-learn}).
The package provides a simple end-to-end workflow going from raw time series data to QoI which can be used for observation-consistent inversion.
Moreover, \verb|LUQ| can also handle noisy prediction data (the process does not change) as we show in one of the numerical examples of Section~\ref{sec:numerics}.

The workflow for using \verb|LUQ| is straightforward. 
We describe a typical use case.
The three initial inputs into the \verb|LUQ| algorithm are an array of $n$ times $\{ t_j\}_{j=1}^n$ (denoted by \verb|times|), a two-dimensional array of $N_{pred}$ ``prediction" time series $\{ \{y^{pred}_{i,j}\}_{j=1}^n\}_{i=1}^{N_{pred}}$ (denoted by \verb|predicted_time_series|), where $y^{pred}_{i,j} = y(t_j)$ for prediction case $i$, and a two-dimensional array of $N_{obs}$ ``observed" time series $\{ \{y^{obs}_{i,j}\}_{j=1}^n\}_{i=1}^{N_{obs}}$ (denoted by \verb|observed_time_series|), where $y^{obs}_{i,j} = y(t_j) + \eta_i(t_{j})$ for observation case $i$.
$\eta_i(t_{j})$ is mean-zero noise, which pollutes the observation data.
These are used to instantiate the LUQ class:
\begin{verbatim}
from luq import LUQ  # Import LUQ module
learn = LUQ(predicted_time_series, observed_time_series, times)
\end{verbatim}

\begin{algorithm}
\SetAlgoLined
\KwData{$m_{min}$, $m_{max}$, $tol$}
\KwResult{$\{\tilde{y}_k\}_{k=1}^{n_{filter}}$}
Let $y_p = \frac{1}{j_f - j_i} \sum_{j = j_i}^{j_f} | y_j |.$ \;
Solve (\ref{eq:spline_opt}) with $m_{min}$ and $m_{min}+1$ to find $S^*_{m_{min}}$ and $S^*_{m_{min}+1}$.\;
Evaluate splines to get $\{y^{m_{min}}_k\}_{k=1}^{n_{filter}}$  and $\{y^{m_{min}+1}_k\}_{k=1}^{n_{filter}}.$ \;
Set $\{y^{old}_k\}_{k=1}^{n_{filter}} = \{y^{m_{min}}_k\}_{k=1}^{n_{filter}}$ and $\{\tilde{y}_k\}_{k=1}^{n_{filter}} = \{y^{m_{min}+1}_k\}_{k=1}^{n_{filter}}$\;
Calculate $Error = \frac{1}{y_p} \sum_{k=1}^{n_{filter}} \left| y^{old}_k - \tilde{y}_k \right| .$\;
Set $m = m_{min} +1.$ \;
\While{$Error > tol$ and $m < m_{max}$}{
Set $m = m +1$ and $y_k^{old} = \tilde{y}_k$  \;
Solve (\ref{eq:spline_opt}) with $m$ to find $S^*_{m_{min}}$. \;
Evaluate spline to get $\{\tilde{y}_k\}_{k=1}^{n_{filter}} = \{y^{m}_k\}_{k=1}^{n_{filter}}$. \;
Calculate $Error = \frac{1}{y_p} \sum_{k=1}^{n_{filter}} \left| y^{old}_k - \tilde{y}_k\right|.$\;
}
\caption{filtering data.}
\label{alg:cleaning}
\end{algorithm}

\subsection{Filtering data}
Next, the data are filtered following Algorithm~\ref{alg:cleaning} for constructing and sampling from a piecewise linear spline with optimally chosen knots.
We first choose a time window over which to filter the data: $[t_{j_i}, t_{j_f}]$, where $1 \leq j_i < j_f \leq n$.
Our goal is to take $n_{filter}$ uniform (in time) filtered data measurements at times $\{ \tilde{t}_k\}_{k=1}^{n_{filter}}$, with  $\tilde{t}_1 = t_{j_i}$ and $\tilde{t}_{n_{filter}} = t_{j_f}$.
Denote the filtered data measurements $\{ \tilde{y}_k\}_{k=1}^{n_{filter}}$.

For each predicted and observed time series, we find optimal piecewise linear splines fitting the raw (possibly noisy) data.
(Note that higher-order splines could, in principal, be used.)
For a spline with $m$ knots, let $\{\bar{t}_k\}_{k=1}^m$ be the locations of the knots and $\mathbf{f} = \{f_k\}_{k=1}^m$ be the values at the knots.
Fix two of the knots at the endpoints of the time window: $\bar{t}_1 = \tilde{t}_1$ and $\bar{t}_m = \tilde{t}_{n_{filter}}$.
Let $\mathbf{\bar{t}} = \{\bar{t}_k\}_{k=2}^{m-1}$, and let $S_m(\mathbf{\bar{t}}, \mathbf{f})$ be the piecewise linear spline defined by those parameters, i.e.,
$$ S_m(\mathbf{\bar{t}}, \mathbf{f}) = \sum_{i=1}^m f_i \ell_i(t), $$
where $\ell_i(x)$, $i=1, \ldots, m$ comprise the linear spline basis:
$$ \ell_i : [\bar{t}_1, \bar{t}_m] \to \mathbb{R}, \quad \ell_i(t) = \begin{cases}
	(t - \bar{t}_{i-1})/(\bar{t}_i - \bar{t}_{i-1}), & t \in [\bar{t}_{i-1}, \bar{t}_i], \\
	(\bar{t}_{i+1} - t)/(\bar{t}_{i+1} - \bar{t}_i), & t \in [\bar{t}_i, \bar{t}_{i+1}], \\
	0, & \text{otherwise}.
\end{cases} $$

Given the number of knots $m$, we find the optimal spline that fits the data in a least squares sense by solving the optimization problem:
\begin{equation}
\mathbf{\bar{t^*}}, \mathbf{f^*} = \argmin_{\mathbf{\bar{t}}, \mathbf{f}} \sum_{j=j_i}^{j_f} \left(y_j - S_m(\mathbf{\bar{t}}, \mathbf{f})(t_j) \right)^2. \label{eq:spline_opt}
\end{equation}
The corresponding optimal spline is denoted $S_m^* := S_m(\mathbf{\bar{t^*}}, \mathbf{f^*})$.
To avoid numerical convergence issues, we add a constraint within to (\ref{eq:spline_opt}) that $\bar{t}_1 \leq \bar{t}_k \leq \bar{t}_m$ for $2 \leq k \leq m$.
In \verb|LUQ|, the constrained optimization problem is solved with a Trust Region Reflective least-squares curve fitting algorithm \cite{Branch:1999:SIC} (the function \verb|curve_fit| from the \verb|scipy.optimize| library in Python).
Let $\{y^m_k\}_{k=1}^{n_{filter}}$ be defined by $y^m_k  = S_m^*(\tilde{t}_k)$.

In order to keep the filter data model as simple as possible, it is preferable to use as few knots as necessary.
Hence, in \verb|LUQ| an adaptive strategy is used to find a simple, yet accurate approximating spline.
Given a minimum number of knots $m_{min}$, a maximum number of knots $m_{max}$, and a tolerance, Algorithm \ref{alg:cleaning} returns the filtered data $\{ \tilde{y}_k\}_{k=1}^{n_{filter}}$.
In Algorithm \ref{alg:cleaning}, optimal splines are formed with successive numbers of knots and are evaluated at the filter times.
When the $1$-norm distance between successive iterations of filtered data (normalized by the absolute average of the raw data in the given time window, which is denoted by $y_p$ in Algorithm~\ref{alg:cleaning}) falls below the specified tolerance or the maximum number of knots is reached, the algorithm terminates, resulting in filtered data: samples from the simplest spline that meets the convergence criteria.
The procedure is performed by iterating over all of the predicted and observed time series data, outputting in the two-dimensional arrays $\{ \{ \tilde{y}^{pred}_{i,k}\}_{k=1}^{n_{filter}} \}_{i=1}^{N_{pred}}$ and $\{ \{ \tilde{y}^{obs}_{i,k}\}_{k=1}^{n_{filter}} \}_{i=1}^{N_{obs}}$.
In other words, optimal splines with possibly different numbers of knots and knot locations are constructed for each distinct time series data within the predicted and observed sets of data as illustrated in Figure~\ref{fig:harmonic-oscillator-cleaning-data}.
Below, we work with these data sets in their matrix forms denoted by $Y^{pred}$ and $Y^{obs}$, respectively.
Within \verb|LUQ|, these matrices along with the vector containing the array of times of the filtered data become attributes of the \verb|LUQ| object. 

This filtering is done within \verb|LUQ| by
\begin{verbatim}
learn.filter_data(time_start_idx=time_start_idx, time_end_idx=time_end_idx,
                 num_filter_obs=num_filter_obs, tol=tol, min_knots=min_knots, 
                 max_knots=max_knots)
\end{verbatim}
where \verb|time_start_idx| is the index of the beginning of the time window, \verb|time_end_idx| is the index of the end of the time window, \verb|num_filter_obs| is the number of uniformly spaced filter observations to take, \verb|tol|, \verb|min_knots|, and \verb|max_knots| are the tolerance, minimum, and maximum number of knots for Algorithm \ref{alg:cleaning}.

\subsection{Clustering and classifying data}
Once the data are filtered, we learn the different types of dynamics present through application of clustering algorithms on the filtered prediction data $Y^{pred}$.
The \verb|scikit-learn| package \cite{scikit-learn} contains several methods for performing clustering on unlabeled datasets that are utilized in \verb|LUQ|.

Perhaps the most popular is the $k$-means algorithm \cite{arthur2006k} which divides data into $k$ clusters of equal variance in a way that minimizes the within-cluster variance.
It is widely adopted because of its simplicity, robustness, and scalability.
The number of clusters, $k$,  must be manually specified.

Another popular choice is to use a Gaussian mixture model (GMM) \cite{mclachlan1988mixture} where the data are assumed to be samples from distribution defined by a weighted mixture of $k$ Gaussian distributions.
The relevant cluster parameters are the means and covariances of each Gaussian along with their corresponding weights.
As with the $k$-means algorithm, the number of clusters must be manually specified.

Spectral clustering \cite{ng2002spectral,von2007tutorial} is at times a useful alternative that uses the spectrum of an affinity matrix between samples to perform dimension reduction.
Following the dimension reduction, a clustering algorithm (e.g., $k$-means or GMM) is then used on the low-dimensional space.
Spectral clustering can be incredibly fast if the affinity matrix is sparse.
It works well for a small number of clusters, but its performance does not scale well for large numbers of clusters.
Again, the number of clusters must be manually specified.

Clustering methods that do not require the number of clusters to be given as input do exist (e.g., Density-Based Spatial Clustering of Applications with Noise (DBSCAN algorithm) \cite{ester1996density, schubert2017dbscan}).
However, they often require that data that are deemed ``noisy" be removed from the data set prior to any clustering.

\verb|LUQ| fully supports the $k$-means, GMM, and spectral clustering routines within \verb|scikit-learn|.
Because of the different methodologies on which each of these clustering algorithms is based, it is difficult to quantitatively compare their outputs.
Hence in \verb|LUQ|, while all are available, the user must choose the type of clustering algorithm to use, along with necessary parameters, e.g., the number of clusters.
Once one of the above clustering routines is used to cluster the filtered prediction data, the data are labeled with an integer designating the cluster to which it belongs. 
Within \verb|LUQ|, the default choice of a clustering algorithm is $k$-means with three clusters and ten random initializations for the cluster centers.

While $k$-means and GMMs implicitly contain a classifier model that can be applied to other data (e.g., the observed data), the other methods discussed do not.
Support Vector Machines (SVMs) are a robust class of supervised learning methods that can be used for  regression, outlier detection, and most notably classification.
SVM classification is effective in high dimensions, memory efficient, and extremely versatile.
A wide range of linear and nonlinear kernel functions can be specified for the decision function, allowing for a wide range of nonlinear separation behavior to be able to be captured.

\verb|scikit-learn| supports a wide variety of SVM classifiers in its \verb|sklearn.svm.SVC| class, which leverages the widely used library \verb|LIBSVM| \cite{chang2011libsvm} for SVMs.
Within the \verb|LUQ| framework, an array of dictionaries of arguments for \verb|sklearn.svm.SVC| are proposed, each defining a class of SVMs (kernels, coefficients, tolerances).
For each proposal, a series of SVM classifiers are trained using the labeled output from the clustering method in order to perform a $k$-fold cross-validation with $k$ series of training and testing sets made from dividing up the samples.
For each proposed SVM class, the $k$-fold cross-validation results in an average misclassification rate (sometimes called classification error rate) defined as the average (over $k$) of the proportion of misclassified test samples.
Whichever proposed SVM class results in the best average misclassification rate is then trained on the entire set of filtered predicted data.
The resulting SVM is subsequently used to label the filtered observed time series data.
Denote the labels as $\{ l^{pred}_i \}_{i=1}^{N_{pred}} $ and $\{ l^{obs}_i \}_{i=1}^{N_{obs}} $, where $l^{pred}_i$ and $l^{obs}_i$ are cluster numbers for predictions and observations $i$, respectively.
The SVM classifier is easily stored, which makes it available to classify any new data sets as they become available.
Within \verb|LUQ| the default proposals for SVM classes are those defined using linear, radial basis function, polynomial, and sigmoid kernels with default coefficients and tolerances and a 10-fold cross-validation; however, the user can propose any possible SVM class available in \verb|sklearn.svm.SVC|.

In \verb|LUQ| the clustering and classification is done by
\begin{verbatim}
learn.dynamics(cluster_method='kmeans',
               kwargs={'n_clusters': 3, 'n_init': 10},
               proposals = ({'kernel': 'linear'},
                {'kernel': 'rbf'}, {'kernel': 'poly'}, {'kernel': 'sigmoid'}),
               k = 10)
\end{verbatim}
where \verb|cluster_method| defines the type of clustering algorithm to use, \verb|kwargs| is a dictionary of arguments for the clustering algorithm, \verb|proposals| is an array of dictionaries of proposed arguments for \verb|sklearn.svm.SVC|, and \verb|k| is the $k$ for the $k$-fold cross-validation.
The printed output gives information about the selection of SVM:
\begin{verbatim}
0.011 misclassification rate for  {'kernel': 'linear'}
0.037 misclassification rate for  {'kernel': 'rbf'}
0.022 misclassification rate for  {'kernel': 'poly'}
0.4170000000000001 misclassification rate for  {'kernel': 'sigmoid'}
Best classifier is  {'kernel': 'linear'}
Misclassification rate is  0.011
\end{verbatim}

\subsection{Feature extraction}
The final step in the \verb|LUQ| framework is feature extraction.
Feature extraction is performed over each cluster of dynamics, rather than the entire data set. 
Let $Y_l^{pred}$ and $Y_l^{obs}$ be submatrices of $Y^{pred}$ and $Y^{obs}$, respectively, for samples labeled in cluster $l$. 
For simplicity assume that $Y_l^{pred}$ and $Y_l^{obs}$ are non-empty for $1 \leq l \leq n_{clusters}$.
In order to avoid numerical errors and improve the quality of the feature extraction, the features are standardized by removing the mean and scaling to unit variance for each $Y_l^{pred}$, and $Y_l^{obs}$ is transformed accordingly.

The \verb|sklearn.decomposition.KernelPCA| class  within \verb|scikit-learn| provides support of a wide variety of kernel PCA methods which are utilized within the \verb|LUQ| framework.
The user specifies an array of proposal arguments defining kernels, coefficients, tolerances, etc., for types of kernel PCA methods supported within this class.
The user also provides either a.) the number of desired QoIs or b.) the minimum proportion of variance that must be explained by the QoIs.
The proportion of variance explained by $n$ QoIs is the ratio of the sum of the eigenvalues associated with the first $n$ principal components with the total some of the eigenvalues.
Each of the proposed kernel PCAs are performed over each $Y_l^{pred}$.
In method a.), for each kernel PCA the proportion of variance associated with the given number of QoIs, $n$, is calculated.
For each cluster $l$, the kernel PCA transform which explains the greatest proportion of variance with $n$ components is selected and applied to both $Y_l^{pred}$ and $Y_l^{obs}$.
The first $n$ components (in each row) of the transformed matrices are the corresponding $n$ QoI.
In method b.), for each kernel PCA, the minimum number of QoI that have a proportion of variance greater than the prescribed minimum is calculated.
The minimum $n$ across the proposals is selected.
If multiple methods result in the same minimum $n$, then the one that explains the most variance is selected.
For each cluster $l$, the chosen kernel PCA is applied to both $Y_l^{pred}$ and $Y_l^{obs}$.
The first $n$ components (in each row) of the transformed matrices are the corresponding $n$ QoI.

Within \verb|LUQ|, the best kernel PCAs are calculated for each cluster and the transformed predictions and observations are computed by
\begin{verbatim}
predict_map, obs_map = learn.learn_qois_and_transform(num_qoi=1,
                             proposals=({'kernel': 'linear'}, {'kernel': 'rbf'},
                             {'kernel': 'sigmoid'}, {'kernel': 'cosine'}))
\end{verbatim}
for case a.), or by
\begin{verbatim}
predict_map, obs_map = learn.learn_qois_and_transform(variance_rate=0.9,
                             proposals=({'kernel': 'linear'}, {'kernel': 'rbf'},
                             {'kernel': 'sigmoid'}, {'kernel': 'cosine'}))
\end{verbatim}
for case b.).
In the above, \verb|num_qoi| is the number of QoIs to use for each cluster, \verb|proposals| is an array of dictionaries of proposed options for \verb|sklearn.decomposition.KernelPCA|, and \verb|variance_rate| is the minimum variance rate that the QoIs need to describe.
The printed output gives information about the selection of the QoIs:
\begin{verbatim}
2 PCs explain 98.9222% of var. for cluster 1 with {'kernel': 'linear'}
2 PCs explain 51.2613% of var. for cluster 1 with {'kernel': 'rbf'}
2 PCs explain 91.4506% of var. for cluster 1 with {'kernel': 'sigmoid'}
2 PCs explain 71.7679% of var. for cluster 1 with {'kernel': 'poly'}
2 PCs explain 98.7201% of var. for cluster 1 with {'kernel': 'cosine'}
---------------------------------------------
Best kPCA for cluster  1  is  {'kernel': 'linear'}
2 PCs explain 98.9222% of variance.
---------------------------------------------
2 PCs explain 99.9406% of var. for cluster 2 with {'kernel': 'linear'}
2 PCs explain 69.4575% of var. for cluster 2 with {'kernel': 'rbf'}
2 PCs explain 94.8345% of var. for cluster 2 with {'kernel': 'sigmoid'}
2 PCs explain 91.5822% of var. for cluster 2 with {'kernel': 'poly'}
2 PCs explain 99.8379% of var. for cluster 2 with {'kernel': 'cosine'}
---------------------------------------------
Best kPCA for cluster  2  is  {'kernel': 'linear'}
2 PCs explain 99.9406% of variance.
---------------------------------------------
\end{verbatim}
\verb|predict_map| and \verb|obs_map| are arrays containing the transformed prediction and observed data for each cluster, i.e., the ``learned" QoIs.
More information about the implementation can be found at \url{https://github.com/CU-Denver-UQ/LUQ} and \cite{mattis_luq}.

%% file: observation-consistent.tex
\section{Applying learned knowledge for observation-consistent inversion} 
\label{sec:inversion}

\subsection{Observation-consistent inversion}

The previous sections describe how the LUQ framework and \verb|LUQ| package transforms each sample of model output data, in the form of a time series, into a QoI associated with each type of dynamical system response.
Ultimately, this creates a mapping from samples of time series of data into samples of QoI for each cluster.
Subsequently, using both the predicted and observed QoI samples, density estimates may be formed on each cluster.
We give a high-level summary of how these density estimates are used to construct an observation-consistent solution to a stochastic inverse problem.
For the interested reader, the works of \cite{BBE11, BES12, BET+14} contain details on the measure-theoretic background for observation-consistent inversion.
For more details on the density-based representation of observation-consistent solutions that we expand upon in this work, we direct the interested reader to \cite{BJW18a}.

Denote by $K$ the number of clusters and $\pspace$ the parameter space.
As illustrated by the plots in the right column of Figure~\ref{fig:harmonic-oscillator-clustering}, the $K$ clusters in data implicitly define a partition of $\pspace$ into $K$ subsets\footnote{In general, the $K$ subsets may be both disconnected and non-convex.}, which we denote by $\pspace_k$ for $k\in\set{1,2,\ldots,K}$.

Although it is not technically necessary to explicitly identify $\pspace_k$ for any $k$, it is computationally trivial to identify which initial set of samples belong to each $\pspace_k$ by simply sorting the labels of the associated sample set of prediction data.
The notation $\pspace_k$ merely provides a formalism that allows us to more easily describe the structure of the solution to the SIP below.
 
For each $k$, let $Q_k(\lambda)$ denote the learned QoI map defined on $\pspace_k$, and let $\mathbb{I}_{\pspace_k}$ denote the indicator function where
\begin{equation}\label{eq:indicator}
	\mathbb{I}_{\pspace_k}(\lambda) = \begin{cases}
											1, & \lambda\in\pspace_k, \\
											0, & \lambda\notin\pspace_k.
										\end{cases}
\end{equation}
Then, for each $k$, let $\predictdenscluster$ denote a predicted density associated with an $L^1$-normalization of $\initdens\mathbb{I}_{\pspace_k}$, where $\initdens$ denotes the initial density used to generate an initial set of parameter samples. 

Similarly, for each $k$, let $\obsdenscluster$ denote the observed density associated with the $k$th cluster of data and $w_k$ denote the {\em observed weights} associated with each cluster so that $0\leq w_k\leq 1$, $\sum_{k=1}^K w_k=1$, and
\begin{equation}
	\sum_{k=1}^K w_k\obsdenscluster(Q_k(\lambda))\mathbb{I}_{\pspace_k}(\lambda)
\end{equation}
defines the observed density over $\pspace$. 

With this notation, we modify the (global) observation-consistent solution described in \cite{BJW18a} as the density given by
\begin{equation}\label{eq:updated-cluster}
	\updatedens(\lambda) = \sum_{k=1}^K w_k\initdens(\lambda)\ratiocluster \mathbb{I}_{\pspace_k}(\lambda),
\end{equation}
where, for each $k$, 
\begin{equation}
	\ratiocluster := \frac{\obsdenscluster(Q_k(\lambda))}{\predictdenscluster(Q_k(\lambda))}.
\end{equation}
Here, we write the observation-consistent solution as an {\em update} to the initial density to make clear the dependence of this solution on how initial samples are generated.
The term $\ratiocluster$ is interpreted as the relative-likelihood that an initial parameter sample $\lambda$ determines a QoI in the $k$th cluster. 
It is useful for performing accept-reject sampling to generate i.i.d.~samples from $\updatedens$. 
It is also useful for constructing a numerical diagnostic on solutions.
Specifically, $\initdens(\lambda)\ratiocluster$ defines an update to the initial density that is re-normalized on each cluster $k$.
In other words, on each cluster $k$, 
\begin{equation}\label{eq:diagnostic}
	1 = \int_\pspace \initdens(\lambda)\ratiocluster = \mathbb{E}_\text{init}(\ratio).
\end{equation}
In this work, we approximate $\predictdenscluster$ and $\obsdenscluster$ using a standard kernel density estimate (KDE) technique with well-established rates of convergence \cite{Terrel_S_JSTOR_1992, Devroye85}.
Thus, computing the sample averages of $\ratiocluster$ and comparing to $1$ provides a numerical diagnostic that these densities are sufficiently accurate on each of the $K$ clusters; see \cite{BJW18a} for more details on this diagnostic.

\subsection{The harmonic oscillator}

\begin{figure}[htbp]
	\centering
	\includegraphics[height=6.5cm]{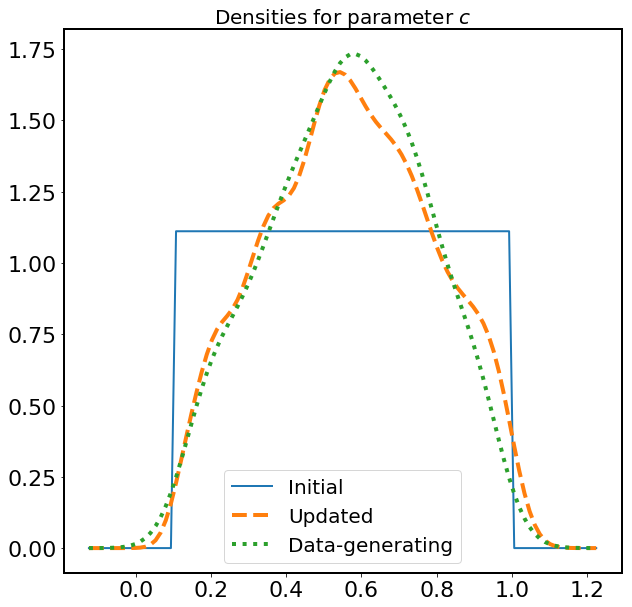}
	\includegraphics[height=6.5cm]{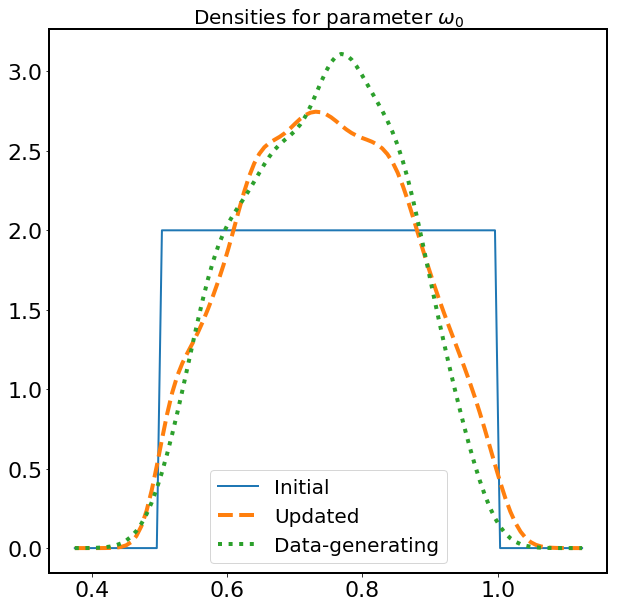}
	\caption{Densities for $c$ (left) and $\omega_0$ (right). The blue solid lines in each plot are the initial uniform densities. The orange dashed lines represent standard weighted kernel density estimates for the updated densities. The green dotted lines are standard kernel density estimates of the data-generating density computed on the finite samples taken from the actual Beta distributions used to construct the density estimates on the learned QoI.}
	\label{fig:harmonic-oscillator-densities}
\end{figure}

For the harmonic oscillator introduced in Section~\ref{sec:LUQ}, a standard KDE is used on each cluster to estimate the predicted and observed densities on the learned QoI.
On clusters 1, 2, and 3, the $\mathbb{E}(\ratiocluster)$ estimates to two significant digits are, 0.98, 0.96, and 0.88,
respectively.
The values for the first two clusters align with what is typically observed in practice for verifying that the density approximations are sufficiently accurate.
The value of 0.88 for the third cluster requires further investigation as it may indicate that either the densities are inaccurately estimated or the predicted density is unable to predict certain observable data.
In this case, the visualizations provided by the classification in data space shown in Figure~\ref{fig:harmonic-oscillator-classifying-data} prove useful. 
It appears that the observed data are within the range of predictions in the third cluster, so we rule out any issues of predictive capacity.
The relatively small number of observed samples that are in the cluster suggests that the most likely cause of this deviation from $1.00$ is the error in the KDE estimate of the observed density.
This is verified by increasing the number of observed samples from $300$ to $1000$ in which case the estimated values of $\mathbb{E}(\ratiocluster)$ on clusters 1, 2, and 3, become 0.95, 1.04, and 1.00, respectively. 
This simply demonstrates the usefulness of the diagnostic in performing a ``sanity check'' on results and in identifying potential sources of errors. 
We continue the analysis below with the relatively low number of $300$ observed samples.

Figure~\ref{fig:harmonic-oscillator-densities} shows several marginal densities for parameters $c$ and $\omega_0$. 
The solid blue curves are the initial uniform densities.
The dashed orange curves show weighted KDEs for the updated densities constructed using~\eqref{eq:updated-cluster}.
Specifically, the KDEs are constructed on the initial samples with weights given by estimates of both $w_k$ and $\ratiocluster$ for $k\in\set{1,2,3}$.
For each $k$, the estimates of $w_k$ are computed using the ratio of number of observed samples classified in cluster $k$ to all observed samples, and $\ratiocluster$ uses the standard KDE estimates obtained for $\obsdenscluster$ and $\predictdenscluster$ evaluated at the number of observed $Q_k$ values. 
Finally, the dotted green curves are standard KDE estimates on the data-generating parameter samples.
We show these KDE estimates of the data-generating distributions to illustrate the impact of finite-sample error in constructing the observed densities.

To better quantify the results, we compute total variation (TV) metrics between densities, i.e., the TV distance between the initial or updated density and the estimated data-generating density for each parameter, as summarized in Table~\ref{tab:harmonic-TV-metrics}.
The TV distance of the updated density estimates from the data-generating densities for each parameter (second column) are reduced by more than 73\% from the distance of the initial density estimates to the data-generating densities (first column).
In other words, the updated density estimates are significantly closer to the data-generating densities. 
Moreover, Table~\ref{tab:harmonic-TV-metrics} also shows the TV distance between KDE approximations of the data-generating density obtained on the finite sample set of $300$ i.i.d. data-generating samples (third column).
Comparing the TV distances in the second and third columns, we conclude that the updated density estimates are comparable to a direct KDE estimate of the data-generating distribution. 


\begin{table}
\centering
\begin{tabular}{|r|r|r|}
\hline
$\| \pi_c^{init} - \pi_c^{DG}  \|_{TV}$  &
$\| \pi_c^{update} - \pi_c^{DG}  \|_{TV} $ & 
$\| \pi_c^{DG} - \pi_c^{DG, exact}  \|_{TV}$ \\
\hline
0.359 & 0.0758 & 0.0872  \\
\hline
\end{tabular}
\vskip 0.2cm
\begin{tabular}{|r|r|r|}
\hline
$\| \pi_{\omega_0}^{init} - \pi_{\omega_0}^{DG}  \|_{TV} $ &
 $ \| \pi_{\omega_0}^{update} - \pi_{\omega_0}^{DG}  \|_{TV} $ &
$\| \pi_{\omega_0}^{DG} - \pi_{\omega_0}^{DG, exact}  \|_{TV} $ \\
\hline
0.372  & 0.100 & 0.0809 \\
\hline
\end{tabular}
\caption{Total variation (TV) metrics for the damped harmonic oscillator problem presented in Section~\ref{sec:LUQ}.
In the first (upper) table, from left-to-right, we report the TV distance between the data-generating marginal density $\pi_c^{DG}$ and (i) the initial marginal density $\pi_c^{init}$, (ii) the updated marginal densities $\pi_c^{update}$, and (iii) the exact marginal distribution $\pi_c^{DG, exact}$, respectively.
In the second (lower) table, from left-to-right, we report similar TV distances for the densities associated with the second parameter, $\omega_0$.}
\label{tab:harmonic-TV-metrics}
\end{table}

\subsection{Impact of different observation times}\label{sec:harmonic-times}

The long-term behaviors of all dynamical responses for the harmonic oscillator involve a steady decay towards equilibrium (i.e., ``eventually'' $y(t)\approx 0$ for sufficiently large $t$).
It is then rather self-evident that the sensitivity of data to parameters decreases over time across all the dynamics.
This is not particularly unique to the harmonic oscillator problem.
It is expected in many dynamical systems that exhibit asymptotic behavior involving equilibrium points or limit cycles.
If data are collected at a time where the dynamics are either no longer sensitive to particular parameter values or only to parameters belonging to certain sets, we do not expect that the QoI extracted from the dynamics  can recover the so-called data-generating distributions.

We assume that a particular useful window of time series data is prescribed, and we focus on the end-to-end analysis of transforming time series data to QoI that provide significant updates of parameter distributions. 
The question of when and where to take measurements in space and time in order to extract useful QoI for observation-consistent inversion is one of optimal experimental design (OED), which is the topic of a future work.
However, in each of the main numerical examples presented in Section~\ref{sec:numerics}, we summarize the impact of using different windows of time data on parameter distribution updates.

%% file: Numerics.tex
\section{Numerical Examples}\label{sec:numerics}


To demonstrate the LUQ framework for solving SIPs arising from dynamical systems, we present three examples of increasing conceptual and computational complexity below.
First, we consider a generalization of the Sel'kov model for glycolysis, a nonlinear system of ordinary differential equations (ODEs), which has a Hopf bifurcation over the parameter domain.
We show the few lines of code required to apply the \verb|LUQ| package, and we also provide the outputs from this package that are reported to the user.
For brevity, we omit the code and output in the second and third examples, but the interested reader may reproduce all of these numerical results using our provided data and Python scripts; see \ref{app:dependencies} for details, including links to our publicly available GitHub and Archive.org repositories.

Next, we turn to the well-known Burgers' equation, a nonlinear partial differential equation (PDE) arising in the study of fluid dynamics, which can produce discontinuities (shock waves) from generic initial conditions.
In this example, uncertainty is assumed in a parameter describing the initial configuration of the wave that subsequently determines how quickly a shock is formed, and results are reported based on temporal measurements taken at different locations in space.
This also demonstrates the overall robustness of the adaptive splines utilized for filtering the data described in Algorithm~\ref{alg:cleaning}.
Specifically, we observe that optimizing the location of an adaptively refined number of knots produces reasonable approximations to discontinuities in time series data.

Finally, we focus on the shallow water equations, a system of nonlinear PDEs arising in the study of coastal circulation and storm surge (coastal flooding).
Here, we consider a physical domain modeling the Shinnecock Inlet located in the Outer Barrier of Long Island, NY, USA.
Both tidal data and a scaling of atmospheric conditions taken from late 2017 through early 2018 are used to simulate an extreme weather event.
The uncertain model inputs describe the parameterization of wind drag, i.e., the parameterization of momentum flux from winds to the water column.
We demonstrate that tidal gauge data measuring water surface elevation can be used to recover key characteristics about the distributions of the uncertain wind drag parameters. 

Within these examples, we also summarize results that motivate on-going and future work involving the proper formulation of ``optimal experimental design'' and ``data assimilation'' problems within this framework.

\input{Example-1}

\input{Example-2}

\input{Example-3}

%% file: Example-1.tex
\subsection{A Hopf bifurcation}\label{sec:Hopf}
In the theory of ODEs, ``Hopf [bifurcations] occur where a periodic orbit is created as the stability of the equilibrium point [changes]'' \cite[pp.~315]{Perko:1991:DED}.
Learning the Quantities of Interest for an ODE with Hopf bifurcations, with its variable dynamics, therefore, is a good demonstration for our proposed framework.

One such class of ODEs arises in the study of biological and biochemical oscillators \cite{Othmer:1977:ECD}, such as the study of the cell cycle \cite{Tyson:1975:CMC} and the study of glycolysis \cite{Selkov:1968:SOG} (the process by which living cells breakdown sugar to obtain energy).
For the purposes of this numerical example, we focus on the following model, which follows from Strogatz \cite[Example 7.3.2]{strogatz2018nonlinear} upon performing a change of variables \cite{munoz2011introduction}:
\begin{align}
x' &= -(x+b) + a \left(y + \frac{b}{a+b^2} \right) + (x+b)^2 \left(y + \frac{b}{a+b^2}\right), \label{eqn:selkov-1} \\
y' &= b-a\left(y+ \frac{b}{a+b^2}\right) - (x+b)^2 \left(y + \frac{b}{a+b^2}\right). \label{eqn:selkov-2}
\end{align}
This model can be framed, for example, as a generalization of the classical Sel'kov model of glycolysis \cite{Selkov:1968:SOG}, where $x$ and $y$ represent concentrations of ADP (adenosine diphosphate) and F6P (fructose 6-phosphate), respectively, and $a,b>0$ are kinetic parameters.
The Hopf bifurcation locus, as a function of $a$, is defined by
\begin{equation} \label{eqn:selkov-locus}
b_1(a) = \sqrt{(1-\sqrt{1-8a}-2a)/2}
\quad \text{and} \quad
b_2(a) = \sqrt{(1+\sqrt{1-8a}-2a)/2}.
\end{equation}
The dynamics of \eqref{eqn:selkov-1}--\eqref{eqn:selkov-2} are classified as follows: If $b<b_1(a)$ or $b > b_2(a)$, then the origin is a stable focus.
If $b_1(a) < b < b_2(a)$, however, then there is a stable periodic orbit.

We are interested in the dynamics of $x$, the concentration of ADP.
The system is solved numerically using the RK45 method \cite{DORMAND198019}, with initial conditions $x(0)=1$ and $y(0)=1$.

For this numerical example, we define data-generating distributions for $a$ and $b$ using two, independent Beta$(2,2)$ distributions over 
$[0.01, 0.124]$ and $[0.05, 1.5]$, respectively to generate a set of 500 samples of time series data from $t=0$ to $t=6.5$ with a measurement rate of $100$ Hz. 
At each observed time, we add a measurement error modeled by an independent and identically distributed $N(0,\sigma^2)$ distribution with $\sigma=0.0125$.
To formulate the predictions at the same measurement frequency, we use independent initial uniform distributions over the parameter intervals. 
We construct 3000 predicted samples using this initial distribution and do not add measurement noise.

\begin{figure}
\centering
\includegraphics[width=0.75\textwidth]{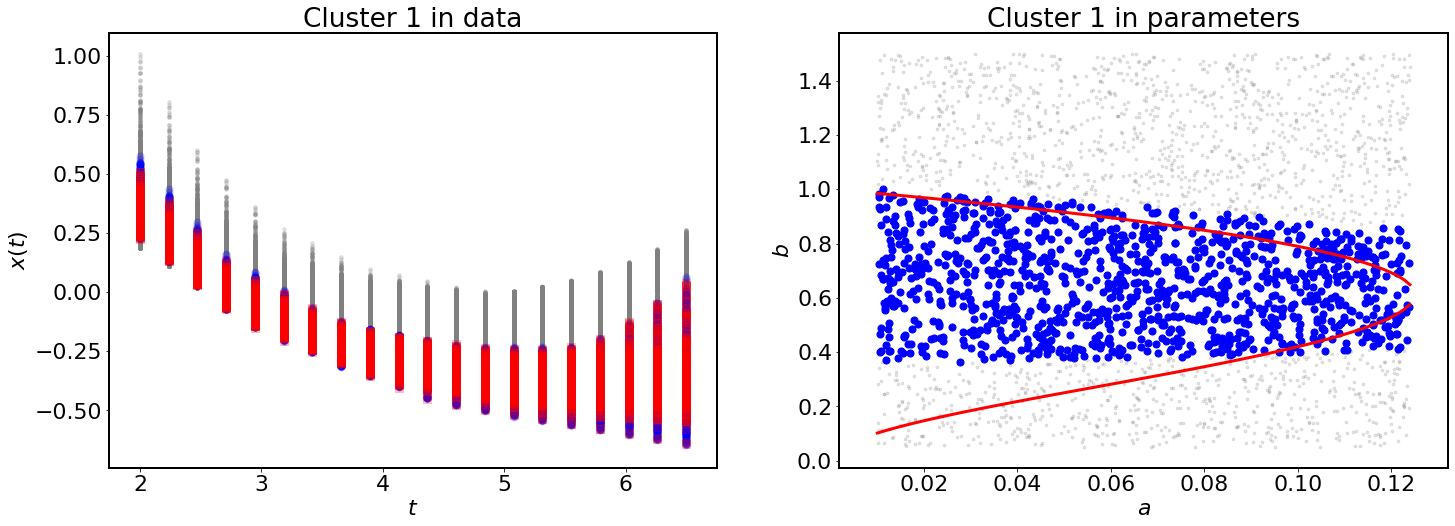}
\includegraphics[width=0.75\textwidth]{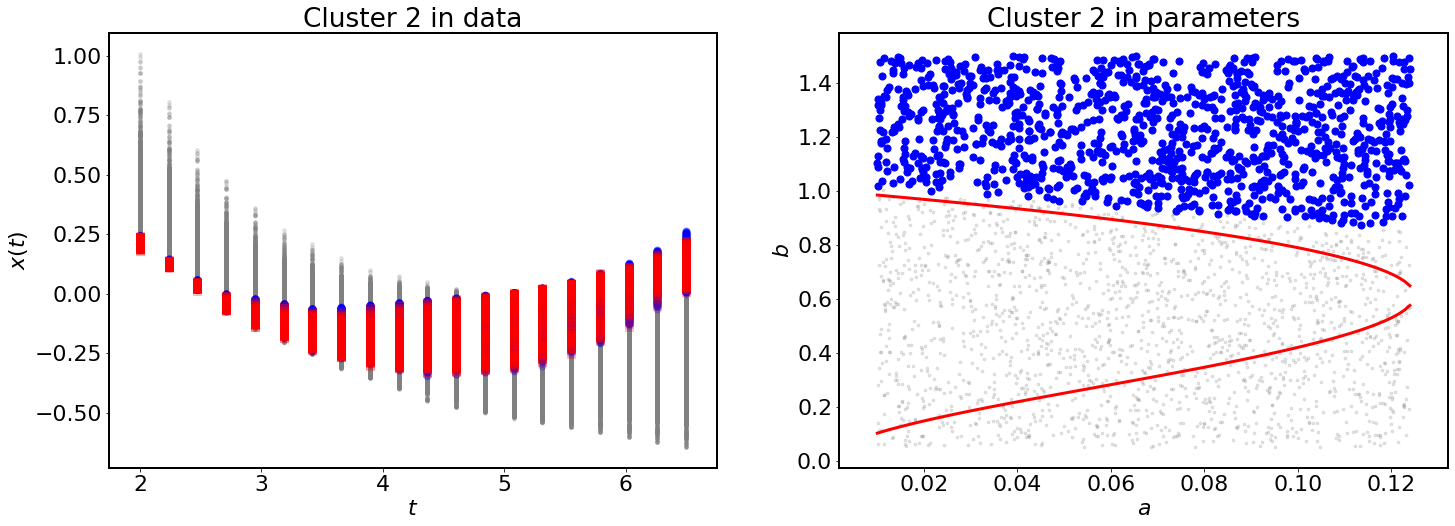}
\includegraphics[width=0.75\textwidth]{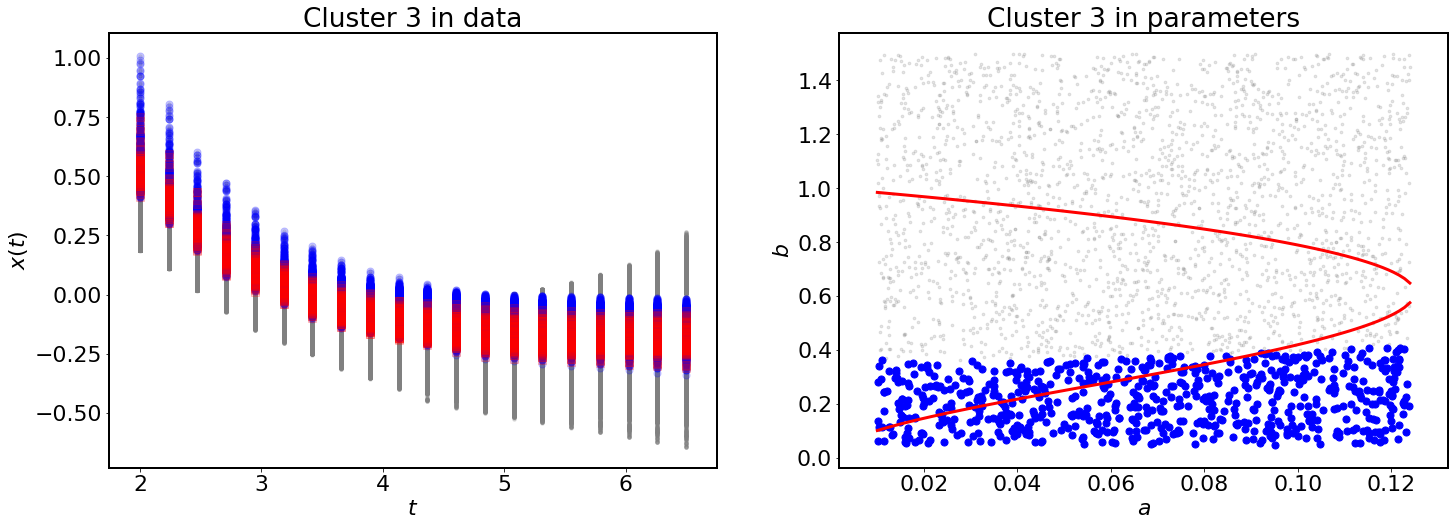}
\caption{Clustering and classification for the Hopf bifurcation problem, Section~\ref{sec:Hopf}.  The large blue dots illustrate the (filtered) predicted data (in time on the left, in the parameter space on the right) using unsupervised learning in the form of $k$-means clustering;  from top to bottom, we show clusters 1, 2, and 3, respectively.  The large red dots (left column) illustrate classified (filtered) observed data.  The smaller gray dots belong to different clusters.  Additionally, the red curves (right column) demonstrate the Hopf bifurcations locus \eqref{eqn:selkov-locus}.}
	\label{fig:hopf-clustering}
\end{figure}

We instantiate a LUQ object, denoted by \verb|LUQ|, and filter the data over a time window of $[2.5, 6.55]$, taking 20 filtered measurements of predicted and observed data using between three and twelve knots:
\begin{verbatim}
learn = LUQ(predicted_time_series, observed_time_series, times)
learn.filter_data(time_start_idx=time_start_idx, time_end_idx=time_end_idx,
                 num_filter_obs=20, tol=5.0e-2, min_knots=3, max_knots=12)
\end{verbatim}
Next, we learn and classify the dynamics using $k$-means clustering with three clusters and the default SVM classifiers:
\begin{verbatim}
learn.dynamics(cluster_method='kmeans', kwargs={'n_clusters': 3, 'n_init': 10})
OUTPUT:
0.003333333333333333 misclassification rate for  {'kernel': 'linear'}
0.009 misclassification rate for  {'kernel': 'rbf'}
0.5083333333333334 misclassification rate for  {'kernel': 'poly'}
0.010666666666666666 misclassification rate for  {'kernel': 'sigmoid'}
Best classifier is  {'kernel': 'linear'}
Misclassification rate is  0.003333333333333333
\end{verbatim}
The resulting SVM has a small average misclassification rate which should result in good accuracy.
The clusters are presented graphically in parameter and data space in Figure \ref{fig:hopf-clustering}.
Observe that three distinct types of dynamics are obtained from the clustering algorithm.
It is obvious that the observed time series (red) are properly being classified by the SVM to match clusters of the predicted time series (blue).
The figures on the right show the predicted samples in parameter space, with the blue dots being the samples in the respective cluster.
The red curves signify the analytically known locations of the Hopf bifurcations.
It is evident that the clustering and classification algorithms are determining regions of parameter space that are predominantly aligned with the known dynamics.

Next, the QoI are learned by choosing kernel PCAs for each cluster of dynamics.
Because there are two uncertain parameters $a$ and $b$, we choose to learn two QoIs.
\begin{verbatim}
predict_map, obs_map = learn.learn_qois_and_transform(num_qoi=2)
OUTPUT:
2 PCs explain 98.9921% of var. for cluster 1 with {'kernel': 'linear'}
2 PCs explain 53.4257% of var. for cluster 1 with {'kernel': 'rbf'}
2 PCs explain 93.9109% of var. for cluster 1 with {'kernel': 'sigmoid'}
2 PCs explain 78.9174% of var. for cluster 1 with {'kernel': 'poly'}
2 PCs explain 97.7487% of var. for cluster 1 with {'kernel': 'cosine'}
---------------------------------------------
Best kPCA for cluster  1  is  {'kernel': 'linear'}
2 PCs explain 98.9921% of variance.
---------------------------------------------
2 PCs explain 95.1797% of var. for cluster 2 with {'kernel': 'linear'}
2 PCs explain 63.3258% of var. for cluster 2 with {'kernel': 'rbf'}
2 PCs explain 95.7184% of var. for cluster 2 with {'kernel': 'sigmoid'}
2 PCs explain 82.1687% of var. for cluster 2 with {'kernel': 'poly'}
2 PCs explain 91.3099% of var. for cluster 2 with {'kernel': 'cosine'}
---------------------------------------------
Best kPCA for cluster  2  is  {'kernel': 'sigmoid'}
2 PCs explain 95.7184% of variance.
---------------------------------------------
2 PCs explain 99.7161% of var. for cluster 3 with {'kernel': 'linear'}
2 PCs explain 57.9535% of var. for cluster 3 with {'kernel': 'rbf'}
2 PCs explain 93.0670% of var. for cluster 3 with {'kernel': 'sigmoid'}
2 PCs explain 77.6821% of var. for cluster 3 with {'kernel': 'poly'}
2 PCs explain 99.3870% of var. for cluster 3 with {'kernel': 'cosine'}
---------------------------------------------
Best kPCA for cluster  3  is  {'kernel': 'linear'}
2 PCs explain 99.7161% of variance.
---------------------------------------------
\end{verbatim}
We see that linear kernels are chosen for clusters 1 and 3 and a sigmoid kernel is chosen for cluster 2.
Moreover, the learned QoI in each cluster all explain a very high proportion of the variance (more than 95\% in each case).
Having learned the QoI for each cluster and transformed the prediction and observed data, we can form weighted KDEs and perform observation-consistent inversion.
A standard KDE is used on each cluster to estimate the predicted and observed densities on the learned QoI.
On clusters 1, 2, and 3, the $\mathbb{E}(\ratiocluster)$ estimates to two significant digits are, 0.99, 0.95, and 1.02,
respectively.

\begin{figure}
\centering
\includegraphics[height=6.5cm]{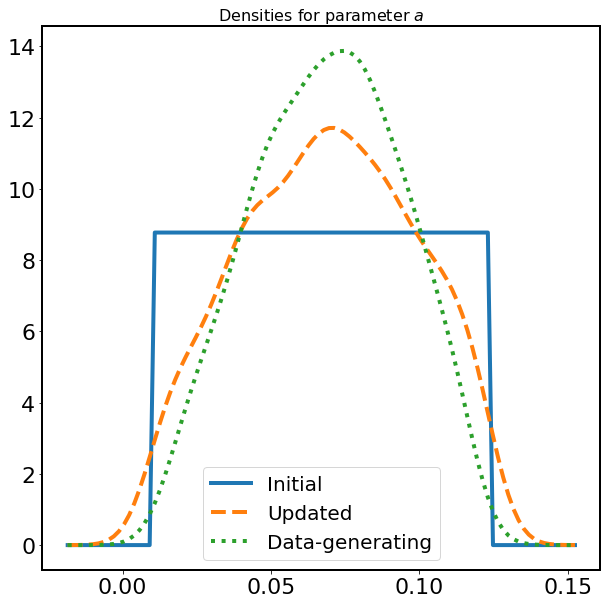}
\includegraphics[height=6.5cm]{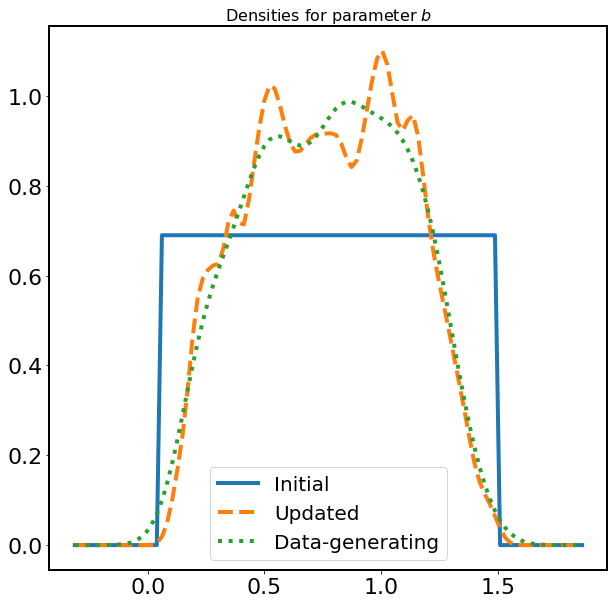}
\caption{Densities for $a$ (left) and $b$ (right). The blue solid lines in each plot are the initial uniform densities. The orange dashed lines represent standard weighted kernel density estimates for the updated densities. The green dotted lines are standard kernel density estimates of the data-generating density computed on the finite samples taken from the actual Beta distributions used to construct the density estimates on the learned QoI.}
	\label{fig:hopf-densities}
\end{figure}

Figure~\ref{fig:hopf-densities} shows several marginal densities for parameters $a$ and $b$. 
The solid blue curves are the initial uniform densities.
The dashed orange curves show weighted KDEs for the updated densities constructed using~\eqref{eq:updated-cluster}.
Specifically, the KDEs are constructed on the initial samples with weights given by estimates of both $w_k$ and $\ratiocluster$ for $k\in\set{1,2,3}$.
For each $k$, the estimates of $w_k$ are computed using the ratio of number of observed samples classified in cluster $k$ to all observed samples and $\ratiocluster$ uses the standard KDE estimates obtained for $\obsdenscluster$ and $\predictdenscluster$ evaluated at the number of observed $Q_k$ values. 
Finally, the dotted green curves are standard KDE estimates on the data-generating parameter samples.
We show these KDE estimates of the data-generating distributions to illustrate the impact of finite-sample error in constructing the observed densities.

To better quantify the results, we compute the TV metrics between densities, i.e., the TV distance between the initial or updated density and the estimated data-generating density for each parameter,
as summarized in Table~\ref{tab:hopf-TV-metrics}.
The TV distance of the updated density estimates from the data-generating densities for each parameter (second column) are reduced by more than 58\% from the distance of the initial density estimates to the data-generating densities (first column).
Moreover, Table~\ref{tab:hopf-TV-metrics} also shows the error in the KDE estimates of the data-generating density resulting from the use of finite-sampling.
We see that the absolute distances between the updated and data-generating density KDEs (second column) are similar to the magnitude of error that arises from using finite-sampling in constructing KDEs (third column).

\begin{table}
\centering
\begin{tabular}{|r|r|r|}
\hline
$\| \pi_a^{init} - \pi_a^{DG}  \|_{TV}$  &
$\| \pi_a^{update} - \pi_a^{DG}  \|_{TV} $ &
$\| \pi_a^{DG} - \pi_a^{DG, exact}  \|_{TV}$ \\
\hline
0.408 &  0.171 & 0.075  \\
\hline
\end{tabular}
\vskip 0.2cm
\begin{tabular}{|r|r|r|}
\hline
$\| \pi_{b}^{init} - \pi_{b}^{DG}  \|_{TV} $ &
$ \| \pi_{b}^{update} - \pi_{b}^{DG}  \|_{TV} $ &
$\| \pi_{b}^{DG} - \pi_{b}^{DG, exact}  \|_{TV} $ \\
\hline
0.322 & 0.077 & 0.060 \\
\hline
\end{tabular}
\caption{Total variation (TV) metrics for the Hopf bifurcation problem, Section~\ref{sec:Hopf}.
In the first (upper) table, from left-to-right, we report the TV distance between the data-generating marginal density $\pi_a^{DG}$ and (i) the initial marginal density $\pi_a^{init}$, (ii) the updated marginal densities $\pi_a^{update}$, and (iii) the exact marginal distribution $\pi_a^{DG, exact}$, respectively.
In the second (lower) table, from left-to-right, we report similar TV distances for the densities associated with the second parameter, $b$.}
\label{tab:hopf-TV-metrics}
\end{table}

We provide some final remarks regarding the time window, $[2.5, 6.55]$, used in generating the above results.
One motivation for using data taken from such an early window of time in the simulation is understood in the context of using this model for predicting spikes in insulin as a response to rising glucose levels and any subsequent medical interventions this may entail. 
From that perspective, it is better to obtain information that is useful for informing parameter values with data taken as early in time as possible. 
Moreover, as time increases, the dynamics of the system tend to either equilibrium or exhibit large spikes depending on the parameters that then appear to return to equilibrium before producing additional spikes. 
Data taken from time windows where all the dynamics produce similar ``leveled'' responses are not as useful for generating QoI sensitive to both parameters. 
For instance, we found that using data from time windows where $t>20$ that are long enough to contain some predicted spikes in glucose levels will produce QoI that are far more sensitive to $b$ than to $a$.
The impact is that updated $b$ parameter distributions are fairly close to the data-generating distribution, but the update to the $a$ parameter distribution is fairly minor by comparison.
Naturally, this provides motivation to see whether we can reformulate our newly proposed framework for the purposes of ``optimal experimental design'', which will be the topic of future work.
However, the interested reader may use our provided data and scripts (see \ref{app:dependencies}) as a starting point to either simulate other data or explore other time windows of data for use in this problem.


%
%
%
%
%
%
%
%
%
%
%
%
%
%
%

%% file: Example-2.tex
\subsection{A shock}\label{sec:Shock}
\begin{figure}
\centering
\includegraphics[width=0.4\textwidth]{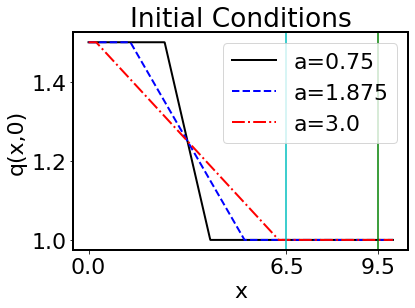}
\caption{Three possible initial conditions $q(x,0)$ for the shock problem, given three choices for the parameter $a$.
Additionally, the measurement locations for the two sets of results are shown in cyan ($x=6.5$) and green ($x=9.5$).}
\label{fig:shock-ic}
\end{figure}
Dynamic partial differential equations (PDEs) are dynamical systems that are often used to model physical laws. 
One such model is the 1D Burgers' equation, a nonlinear PDE arising in the study of fluid dynamics:
$$q_t + \frac{1}{2} (q^2)_x = 0.$$
The spatial domain is the interval $[0, 10]$.
We impose an initial condition of the form
\begin{equation*}
q(x,0) = \begin{cases} 
      f_l, & 0 \leq x\leq 3.25 -a,  \\
       \frac{1}{2} ((f_l + f_r) - (f_l - f_r) \frac{(x-3.25)}{a}), & 3.25 -a < x \leq 3.25 + a, \\
      f_r, & 3.25 + a < x \leq 10,
   \end{cases}
\end{equation*}
where $a \in [0.75, 3]$ is an uncertain parameter and $f_l$ and $f_r$ are positive constants with $f_l > f_r$. 
In this example, we take $f_l = 1.5$ and $f_r = 1$; see, for example, Figure~\ref{fig:shock-ic}.
We assume non-reflecting boundary conditions, allowing waves to pass out of the boundaries without reflection.

This system often can develop discontinuous solutions (shock waves), which complicates calculating a numerical solution. 
We use PyClaw \cite{leveque1997wave,ketcheson2012pyclaw} to calculate weak solutions to the system using a Godunov-type finite volume method with an appropriate limiter and Riemann solver. 
We use a uniform mesh with 500 elements.
The system described above forms a shock at $t = \frac{2a}{f_l - f_r}$.
The shock speed is $\frac{1}{2}(f_l + f_r)$.
We calculate the time series solution at $x=6.5$, i.e., $q(7,t)$ and $x=9.5$ at 1000 evenly spaced time steps between $t=0$ and $t=10$.
We assume a data-generating distribution of $a$ defined by a  Beta$(2,2)$ distribution over 
$[0.75, 3]$ to generate a set of 500 samples of time series data from $t=0$ to $t=10$ with a measurement rate of $100$ Hz. 
At each observed time, we add measurement error that is independent and identically distributed according to an $N(0,\sigma^2)$ distribution with $\sigma=0.025$.
An initial uniform distribution is assumed over the same interval and with the same measurement frequency. 
We take 1000 predicted samples from this distribution for which we also add the same amount of measurement noise as assumed in the predictions.

\begin{figure}
\centering
\includegraphics[width=0.75\textwidth]{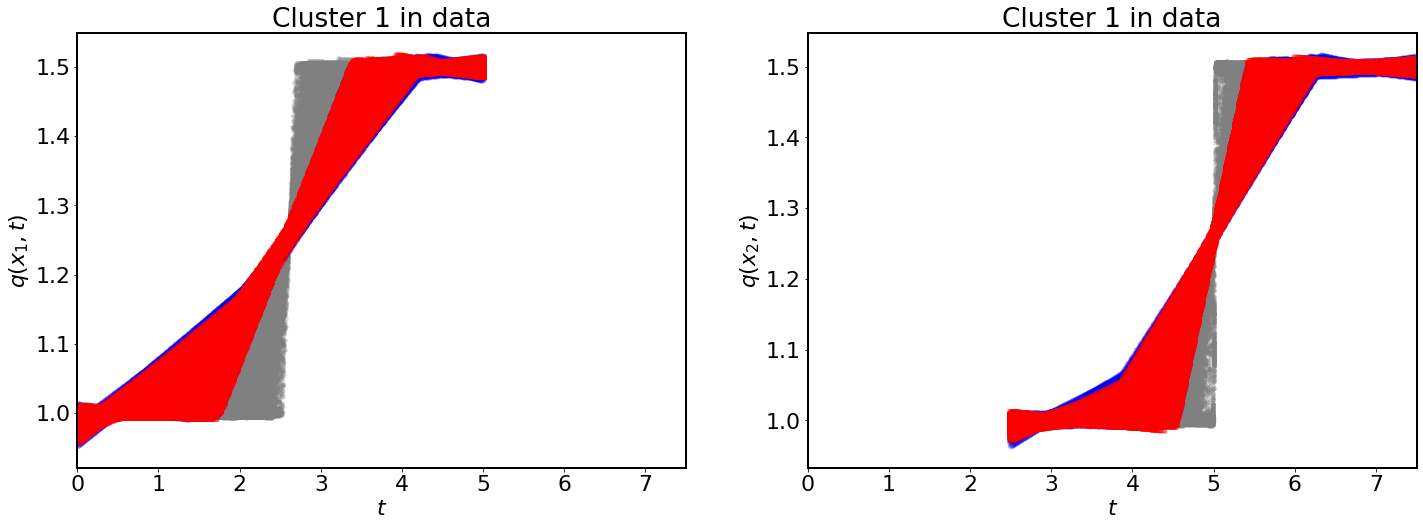}
\includegraphics[width=0.75\textwidth]{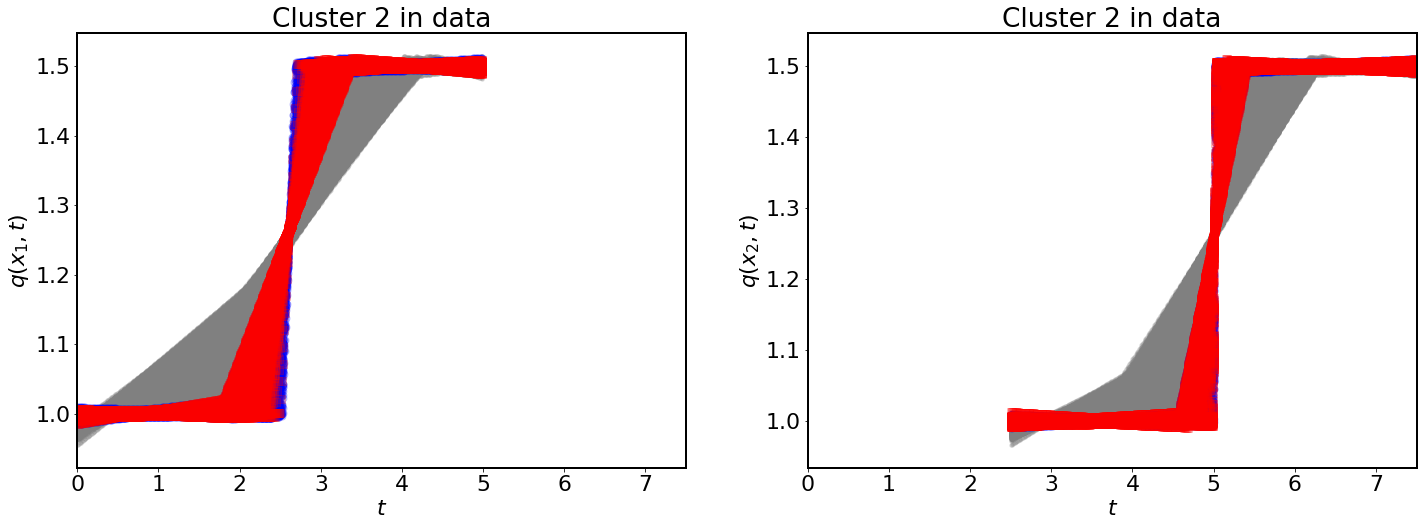}
\caption{Clustering and classification for the shock problem, Section~\ref{sec:Shock}: First column shows results for the first experiment with measurement location $x=6.5$, while for the second column $x=9.5$. The large blue dots (largely obscured by red) illustrate the clustered (filtered) predicted data; from top to bottom, we show clusters 1 and 2, respectively. Moreover, the large red dots illustrate the classified (filtered) observed data, while the smaller gray dots belong to the other cluster.}
\label{fig:shock-clustering}
\end{figure}

First, we analyze the data associate with measurements taken at $x=6.5$. 
We instantiate a LUQ object \verb|LUQ| and filter the data over a time window of $[0, 5]$, taking 500 filtered measurements of predicted and observed data using between three and ten knots.
We expect two main types of dynamics, corresponding to whether or not the wave has transformed into a shock wave by the time it reaches $x=6.5$.
Therefore, we learn and classify the dynamics using $k$-means clustering with two clusters and the default SVM classifiers.
\verb|LUQ| chooses a linear kernel SVM which causes a misclassification rate of $.011$.
Figure \ref{fig:shock-clustering} (left column) shows the two classification clusters. 
Cluster 1 shows dynamics that are far away from forming shocks, and cluster 2 shows dynamics of waves that are almost shocks.
In fact, at $x=6.5$, none of the initial waves defined by the range of $a$ values used here actually form shocks by the time the wave passes this location.

\begin{figure}
\centering
\includegraphics[width=0.75\textwidth]{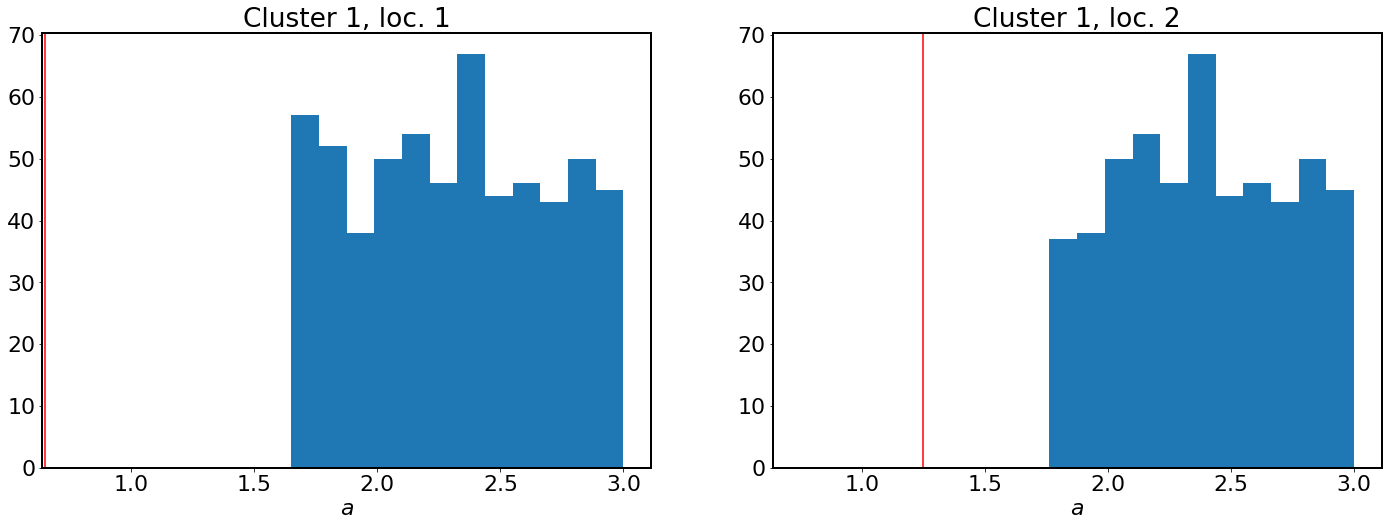}
\includegraphics[width=0.75\textwidth]{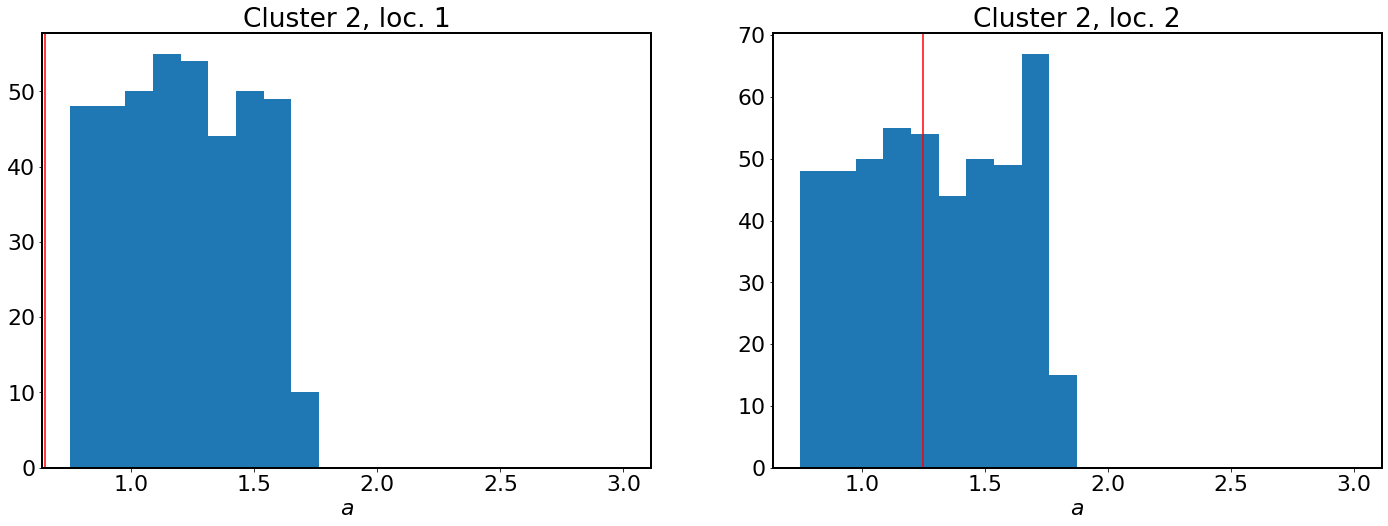}
\caption{Histograms of the parameter $a$ inferred from the clustering of predicted (filtered) data for measurements taken at $x=6.5$ (left column) versus $x=9.5$ (right column). The red line denotes that values to the left have formed a shock.}
\label{fig:shock-params}
\end{figure}

Figure \ref{fig:shock-params} (left column) shows the predicted samples in parameter space.
We see, as expected, smaller values of $a$ in one cluster (cluster 2) and larger values in the other (cluster 1):
this matches what one might assume from the physics.
Next, the QoI are learned by choosing kernel PCAs for each cluster of dynamics.
Because there is only one parameter, $a$, we choose to learn one QoI.
\verb|LUQ| chooses a sigmoid kernel for cluster 1, which explains $67.5589\%$ of the variance and a linear kernel for cluster 2, which explains $29.2163\%$ of the variance.
A standard KDE is used to estimate the predicted and observed densities on the learned QoI for each cluster.
On clusters 1 and 2, the $\mathbb{E}(\ratiocluster)$ estimates to two significant digits are,
1.0 and 0.99,
respectively.

\begin{figure}
\centering
\includegraphics[width=0.4\textwidth]{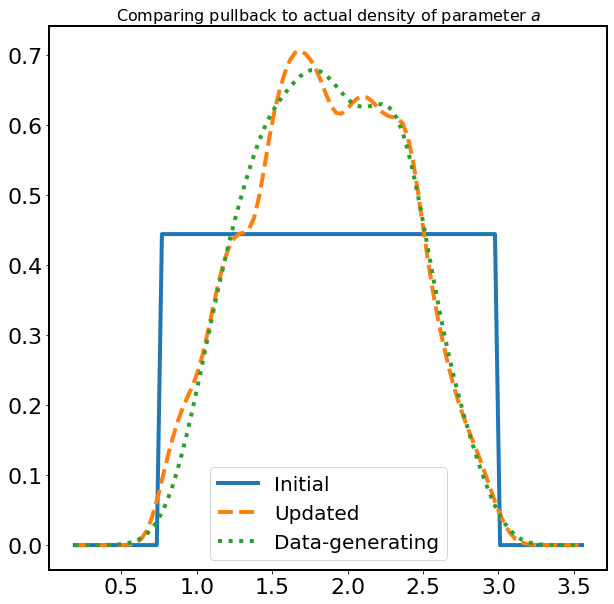}
\includegraphics[width=0.4\textwidth]{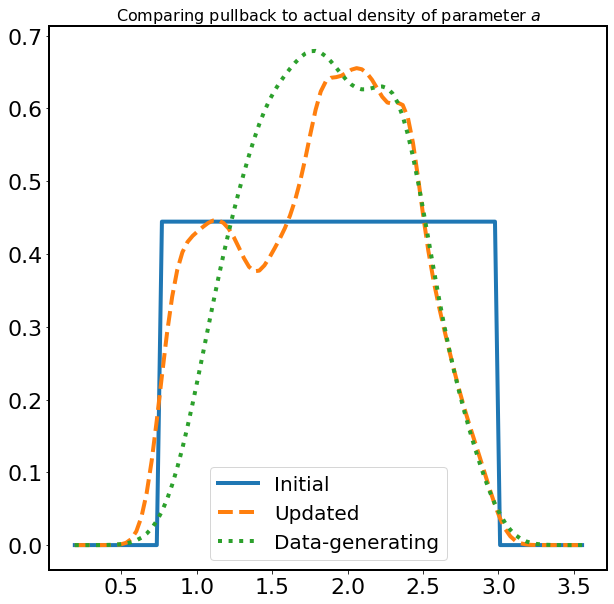}
\caption{Densities for $a$ with data collected at $x=6.5$ (left) and $x=9.5$ (right). The blue solid lines in each plot are the initial uniform densities. The orange dashed lines represent standard weighted kernel density estimates for the updated densities. The green dotted lines are standard kernel density estimates of the data-generating density computed on the finite samples taken from the actual Beta distributions used to construct the density estimates on the learned QoI.}
	\label{fig:shock-densities}
\end{figure}

The plots of Figure~\ref{fig:shock-densities} (left) show several probability densities over $a$. 
In the left plot (for data obtained at $x=6.5$), the updated density matches very well with the data-generating density.
We further analyze the accuracy by looking at the total variation in Table \ref{tab:shock-TV-metrics}.
We see an $87\%$ reduction in the total variation distance of the updated to data-generating densities (second column) compared to the distance of the initial to data-generating densities (first column).
Moreover, the total variation distance of the updated density to the data-generating density (second column) is at the same level as a distance of the direct KDE approximation of the data-generating density (fourth column).
Hence, the updated density is almost as accurate as is possible.

We again emphasize that the results discussed above are for the data collected at measurement location $x=6.5$, where none of the samples from the initial density have actually formed shocks.
To understand the impact of shocks on the data and subsequently on the updated density for $a$, we instead look at measurements taken at $x=9.5$. 
At this location, many of the samples from the initial density have formed shocks.
Once a shock has formed, since the shock speeds are identical for each sample, it is impossible to identify the precise value for the parameter $a$ other than it belonging to a particular set of plausible values.
In other words, any value of $a$ that has caused a shock wave to form produces an identical (up to noise) time series.
The impact of this is analyzed below.

We instantiate a LUQ object \verb|LUQ| and filter the data over a time window of $[2.5, 7.5]$, taking 500 filtered measurements of predicted and observed data using between three and ten knots.
We learn and classify the dynamics using $k$-means clustering with two clusters and the default SVM classifiers.
\verb|LUQ| chooses a linear kernel SVM which causes a misclassification rate of $.001$.
Figure \ref{fig:shock-clustering} (right column) shows the two classification clusters. 
Cluster 1 shows dynamics that are far away from forming shocks, and cluster 2 shows dynamics of waves that are either shocks or almost shocks.
Figure \ref{fig:shock-params} (right column) shows the predicted samples in parameter space.
We see, as expected, in cluster 1 the values of $a$ are generally above the threshold of forming shocks ($a=1.25$), and in cluster 2 the values are below or only slightly above the threshold.
Next, the QoI are learned by choosing kernel PCAs for each cluster of dynamics.
Because there is only one parameter, we again choose to learn one QoI.
\verb|LUQ| chooses a sigmoid kernel for cluster 1, which explains $51.0022\%$ of the variance and a linear kernel for cluster 2, which explains $28.3004\%$ of the variance.
A standard KDE is used on each cluster to estimate the predicted and observed densities on the learned QoI.
On clusters 1 and 2, the $\mathbb{E}(\ratiocluster)$ estimates to two significant digits are
1.0 and 0.98,
respectively.

The right plot in Figure~\ref{fig:shock-densities} shows several probability densities over $a$ based on measurements at $x=9.5$.
The updated density matches very well with the data-generating density for the region where shocks have not formed given by $a>1.25$.
However, for $a\leq 1.25$, it does not match well for the reasons discussed above. 
In Table~\ref{tab:shock-TV-metrics}, we observe the total variation distance of the subsequent updated distribution to the data-generating distribution (third column) is more than when using data from $x=6.5$ (second column), which is consistent with the above analysis of the results.
However, the probability with respect to the updated distribution of the event encompassing cluster 2 is $0.446$ and the probability of the data-generating distribution is 0.448, and error of less than $1\%$.
Hence, the probability of the shock and near-shock event has been computed very accurately.
In other words, while we cannot expect the updated density to accurately describe differences in relative likelihoods in the event defined by $a\leq 1.25$, it can still be used to accurately compute the probability of this event. 

\begin{table}
\centering
\begin{tabular}{|r|r|r|r|}
\hline
$\| \pi^{init} - \pi^{DG}  \|_{TV}$ &  
$\| \pi^{update, 6.5} - \pi^{DG}  \|_{TV} $ & 
$\| \pi^{update, 9.5} - \pi^{DG}  \|_{TV} $ &
$\| \pi^{DG} - \pi^{DG, exact}  \|_{TV}$ \\ 
\hline
0.430 & 0.054 & 0.192 & 0.059 \\
\hline
\end{tabular}
\caption{Total variation (TV) metrics for the shock problem, Section~\ref{sec:Shock}. 
From left-to-right, we report the TV distance between the data-generating density $\pi^{DG}$ and (i) the initial density $\pi^{init}$, (ii) the updated density $\pi^{update,6.5}$ for the first experiment with data collected at $x=6.5$, (iii) the updated density $\pi^{update,9.5}$ for the second experiment ($x=9.5$), and (iv) the exact marginal distribution $\pi^{DG, exact}$, respectively.}
\label{tab:shock-TV-metrics}
\end{table}

%
%
%
%
%
%
%
%
%
%
%
%
%
%
%

%% file: Example-3.tex
\subsection{Wind drag}\label{sec:drag}

In coastal circulation and flooding modeling, a common fluid dynamics approximation is the {\em shallow water approximation}.
If the horizontal length scale of motion (e.g., wavelength) is much greater than the height of the water column, then, as first studied by A.J.C. Barr{\'e} de Saint-Venant in a 1D setting \cite{Saint-Venant:1871:SWE}, one can derive a system of partial differential equations for the evolution of water surface elevation and depth-averaged momentum (in contrast with velocity and pressure as in the Navier--Stokes equations), called the Shallow Water Equations (SWE).
See \cite{Stoker:1957:WWM,Vreugdenhil:1994:NMS,LeVeque:2002:FVM} for a selection of derivations.

With meteorological forcing (wind speed and air pressure), the SWE can be written as
\begin{align}
 	\frac{\partial\zeta}{\partial t} + \nabla \cdot (\mathbf{U} H) &= 0, \label{eqn:wd-swe-1} \\
 	\frac{\partial\mathbf{U}}{\partial t} + \mathbf{U}\cdot\nabla\mathbf{U} + f\mathbf{k}\times\mathbf{U} &= -\nabla\left(\frac{p_s}{\rho_0} + g\zeta\right) + \frac{\boldsymbol\tau_s - \boldsymbol\tau_b}{\rho_0 H}, \label{eqn:wd-swe-2}
\end{align}
for the unknown free surface elevation  $\zeta = \zeta(x, y, t) \in \mathbb{R}$ (alt., water elevation or sea level) and depth-averaged velocity vector $\mathbf{U} = \mathbf{U}(x,y,t) \in \mathbb{R}^2$.
In \eqref{eqn:wd-swe-1}--\eqref{eqn:wd-swe-2}, $H = h + \zeta$ is the total water depth, where $h$ is the still water depth;
$f = 2\Omega\sin\phi$ is the Coriolis parameter, where $\Omega$ is the angular speed of the Earth and $\phi$ is latitude;
$p_s$ is the atmospheric pressure at the free surface and $\rho_0$ the reference density of water;
$g$ is acceleration due to gravity;
$\boldsymbol\tau_s = \rho_a C_d \mathbf{u} ||\mathbf{u}||$ and $\boldsymbol\tau_b$ are the free surface and seabed stresses, respectively;
$C_d$ is the wind drag;
and $\mathbf{u}$ is the wind speed at 10-m.
(In this work, a hybrid friction law is used for the seabed stress, see \cite{ADCIRC:v53}.)

The ADvanced CIRCulation (ADCIRC) coastal ocean model is a continuous-Galerkin, finite-element model of the SWE \cite{Luettich:1992,Westerink:1992}, in which the Generalized Wave Continuity Equation \cite{Lynch:1979} (an equivalent formulation of the SWE) is discretized in space using piecewise-linear elements on unstructured (triangular) grids.
It is used quasi-operationally for coastal engineering applications such as hurricane storm surge hindcasting \cite{Bunya:2010:HRC,Dietrich:2010:HRC,Dietrich:2011:HGW} and forecasting \cite{Dietrich:2013:RTF} and uncertainty quantification \cite{Butler:2015:DSS,Graham:2015:AMT,Graham:2017:MTA}, and can run in both single core and distributed computing environments \cite{Tanaka:2010,Dietrich:2012:PUM}.

Significant uncertainty in these applications exists within the meteorological components, such as the forecasted wind and pressure fields, and the wind drag.
A common formulation of wind drag in storm surge applications is
$C_d = \min[10^{-3}(0.75 + .067 ||\mathbf{u}||), .0025]$,
where the exact numerical coefficients may vary between applications.
Briefly, this formulation for $C_d$ models the observed linear increase with low wind speeds and the ``cut-off'' or ``saturation'' for high wind speeds, see \cite{Garratt:1977:RDC,Letchford:2009:WWS}.

As a numerical experiment for LUQ and observation-consistent inversion, we propose a generalization of the above equation, namely,
$$ C_d = \min[10^{-3}(0.75 + \lambda_1 ||\mathbf{u}||), \lambda_2]. $$
We suppose, for the purposes of the experiment, that the uncertain parameters $(\lambda_1, \lambda_2)$ might lie within $\pm50\%$ of the aforementioned coefficients $.067$ and $.0025$, respectively:
\begin{equation}
	\lambda_1 \in \Lambda_1 := [.0335, .1005], \quad \lambda_2 \in \Lambda_2 := [.00125, .00375]. \label{eqn:wd-lambda}
\end{equation}
Our goal will be to recover an approximation of the distribution on $(\lambda_1, \lambda_2)$ from time series of water surface elevation $\zeta(x_0, y_0, t)$ at a given, fixed location $(x_0, y_0)$.

\begin{figure}
\centering
\includegraphics[height=7cm]{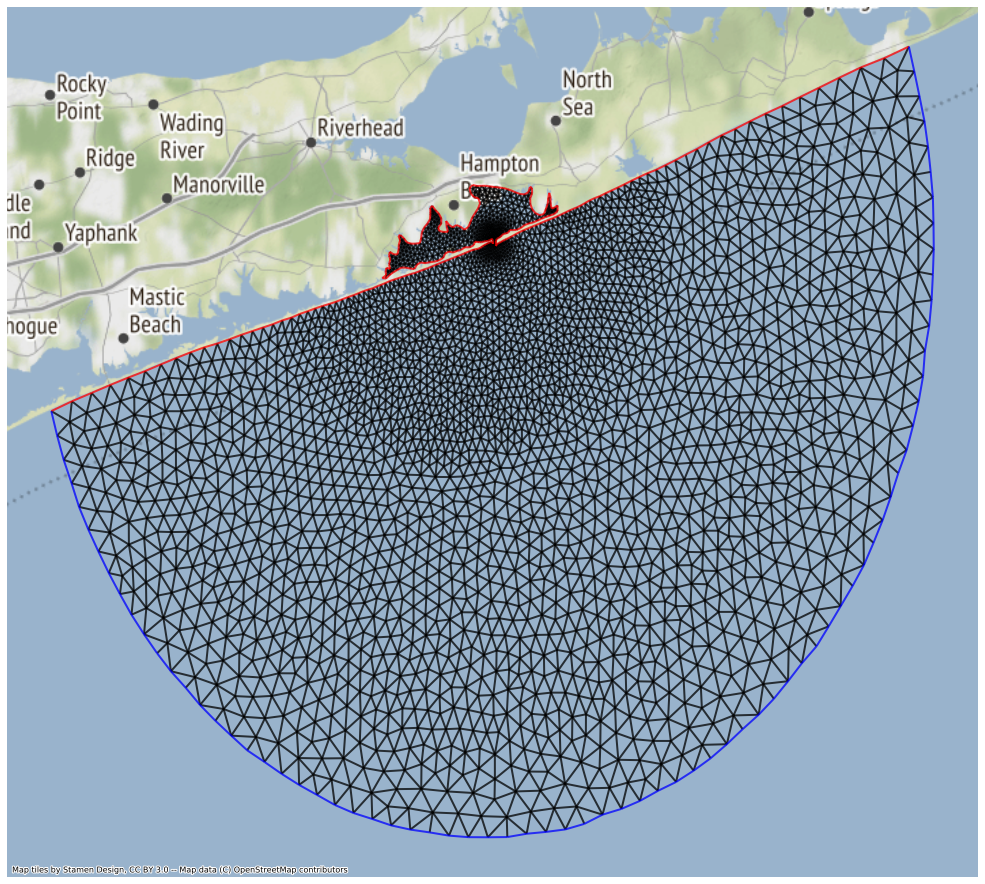}%
\raisebox{3.5cm}{\includegraphics[height=3.5cm]{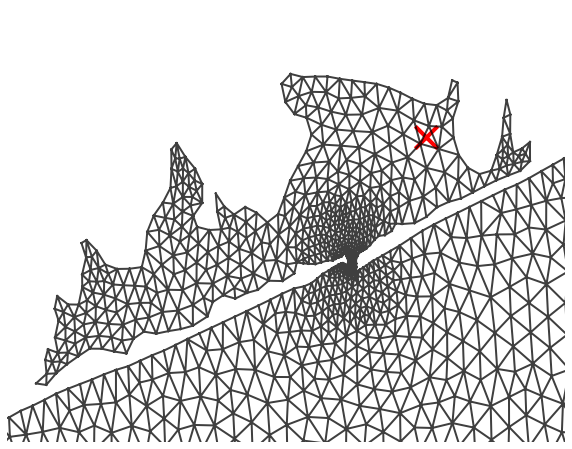}}
\caption{Mesh for the Shinnecock Inlet, with tidal forcing boundary (blue), no-flow boundary (red), and the time series measurement location $(x_0, y_0)$ (see inset, red ``$\times$''). The map tiles, by Stamen Design, are licensed under CC BY 3.0. Map data, \textcopyright~OpenStreetMap contributors, are licensed under ODbL.}
\label{fig:wd-si-triplot}
\end{figure}

We configure an ADCIRC model using the well-tested Shinnecock Inlet mesh \cite{ADCIRC:2020:SI} (modeling an inlet in the Outer Barrier of Long Island, NY, USA), with external forcing given by tides, winds, and constant air pressure for a period of 34 days (29 December 2017 -- 31 January 2018).
The relatively small mesh (approx. 3.1k nodes and 5.8k triangles) is shown in Figure~\ref{fig:wd-si-triplot}.
The tidal forcing is reconstructed from the TPXO9.1 harmonic tidal constituents \cite{Egbert:2002:EIM} using OceanMesh2D \cite{Roberts:2019:OM2D}, the air pressure is assumed spatially and temporally constant (1013 millibars), and the free surface stress from winds is computed from 0.25$^\circ$ hourly CFSv2 10-m wind fields \cite{Saha:2014:CFS} interpolated onto the mesh.
The winds are modified for the purposes of our numerical experiment.
First, to simulate a more extreme event (winds up to Category 4), they are artificially scaled by a factor of three.
Second, to reduce spurious, non-physical oscillations near the outer boundary, they are smoothly reduced to zero near the outer boundary.

For the numerical experiment, we use two, independent data-generating distributions, $\lambda_1 \sim \text{Beta}(5,2)$ and $\lambda_2 \sim \text{Beta}(1.5,7)$ (translated and scaled to the intervals $\Lambda_1$ and $\Lambda_2$, respectively), to generate 100 samples of time series data (water surface elevation) measured every minute for the month of January 2018 at $(x_0, y_0) = (-72.45^\circ, 40.87^\circ)$, see Figure~\ref{fig:wd-si-triplot}.
To each observation, we add independent, normally distributed measurement error $N(0,\sigma^2)$ with $\sigma = .005$.
We additionally use uniform distributions on identical intervals, $\lambda_1 \sim \text{unif}(\Lambda_1)$ and $\lambda_2 \sim \text{unif}(\Lambda_2)$, to generate 1000 predicted (noise-free) time series.

First, we instantiate a LUQ object \texttt{LUQ} and filter the data over a time window of [03 Jan 2018 18:41:00 GMT, 03 Jan 2018 22:00:00 GMT], represented in the data using Unix time (seconds since 1 January 1970 0:00) $T_1 = [1515004860, 1515016800]$.
We take ten filtered measurements of predicted and observed data using between five and ten knots.
We do not expect significantly different dynamics based on the choice of parameters $(\lambda_1, \lambda_2)$, so we do not perform clustering, do not construct an SVM classifier, etc.
Because there are two uncertain parameters, we choose to learn two QoI.
LUQ chooses a linear kernel PCA, which explains approximately $89\%$ of the variance in the filtered data.
A standard KDE is used to estimate the predicated and observed densities on the learned QoI.
The $\mathbb{E}(\ratio)$ estimate to two significant digits is $0.83$.

\begin{figure}
	\centering
	\includegraphics[height=6.5cm]{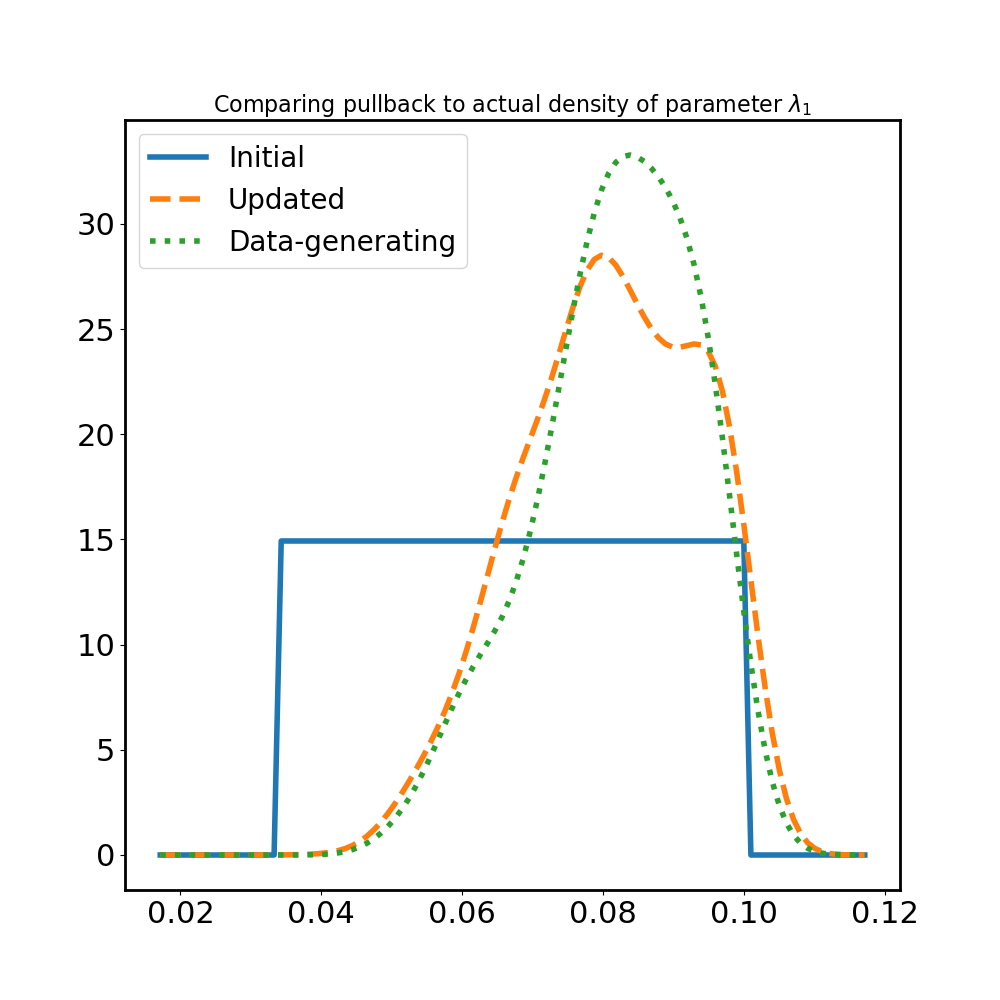}
	\includegraphics[height=6.5cm]{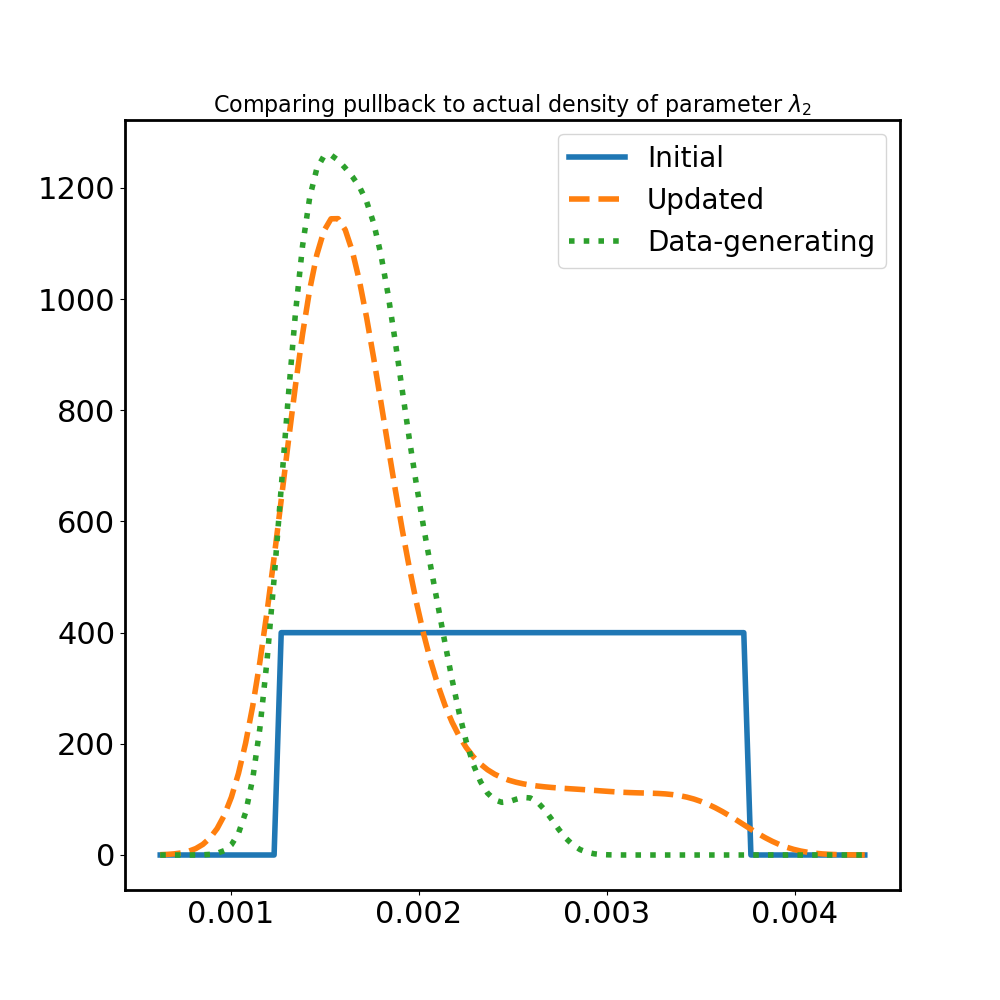}
	\caption{Densities for $\lambda_1$ (left) and $\lambda_2$ (right) for the first numerical experiment presented in Section~\ref{sec:drag}.
	The blue solid lines in each plot are the initial uniform densities.
	The orange dashed lines represent standard weighted kernel density estimates for the updated densities.
	The green dotted lines are standard kernel density estimates of the data-generating density computed on the finite samples taken from the actual Beta distributions used to construct the density estimates on the learned QoI.}
	\label{fig:wd-densities-v1}
\end{figure}

Figure~\ref{fig:wd-densities-v1} shows marginal probability densities over $\lambda_1$ (left) and $\lambda_2$ (right).
The solid blue curves are the initial uniform densities.
The dashed orange curves show weighted KDEs for the updated densities with $K=1$ and $w_k=1$.
Specifically, the KDEs are constructed on the initial samples with weights given by estimates of $\ratio$, using the standard KDE estimates obtained for $\obsdens$ and $\predictdens$ evaluated at the number of observed $Q$ values.
Finally, the dotted green curves are standard KDE estimates on the data-generating parameter samples.

\begin{table}
\centering
\begin{tabular}{|r|r|r|}
\hline
$\| \pi_{\lambda_1}^{init} - \pi_{\lambda_1}^{DG}  \|_{TV}$  &
$\| \pi_{\lambda_1}^{update} - \pi_{\lambda_1}^{DG}  \|_{TV} $ & 
$\| \pi_{\lambda_1}^{DG} - \pi_{\lambda_1}^{DG, exact}  \|_{TV}$ \\
\hline
0.744 & 0.165 & 0.093 \\
\hline
\end{tabular}
\vskip 0.2cm
\begin{tabular}{|r|r|r|}
\hline
$\| \pi_{\lambda_2}^{init} - \pi_{\lambda_2}^{DG}  \|_{TV} $ &
 $ \| \pi_{\lambda_2}^{update} - \pi_{\lambda_2}^{DG}  \|_{TV} $ &
$\| \pi_{\lambda_2}^{DG} - \pi_{\lambda_2}^{DG, exact}  \|_{TV} $ \\
\hline
1.061 & 0.276 & 0.160 \\
\hline
\end{tabular}
\caption{Total variation (TV) metrics for the first numerical experiment presented in Section~\ref{sec:drag}, using the time window $T_1$.
From left-to-right, in each of the two tables, we report the TV distances between the data-generating marginal densities $\pi_{\lambda_i}^{DG}$ ($i=1$ in the first table and $i=2$ in the second) and (i) the initial marginal densities $\pi_{\lambda_i}^{init}$, (ii) the updated marginal densities $\pi_{\lambda_i}^{update}$, and (iii) the exact marginal distribution $\pi_{\lambda_i}^{DG, exact}$ , respectively.}
\label{tab:wd-TV-metrics-1}
\end{table}

The updated densities match very well with the data-generating densities.
We further analyze the accuracy by looking at the total variation in Table~\ref{tab:wd-TV-metrics-1}.
We see a reduction of more than $73\%$ in the total variation distance of the updated to data-generating densities (second column) compared with the distance of the initial to data-generating densities (first column) for each parameter.
Moreover, the total variation distances of the updated marginal densities to the data-generating marginal densities (second column) are similar to the distances obtained from a direct KDE approximation of the data-generating marginal densities (third column).
Taking into account the aforementioned $\mathbb{E}(\ratio)$ estimate together with the measured TV reduction, we can confidently say that the stochastic inverse problem (SIP) has been solved accurately for both $(\lambda_1, \lambda_2)$.

To show the impact of experimental design, we conduct a second numerical experiment with the same data set but a new time window.
We instantiate a LUQ object and filter data over a time window [02 Jan 2018 00:00:00 GMT, 02 Jan 2018 23:59:00 GMT], represented in the data using Unix time $T_2 = [1514851200, 1514937540]$.
In this time window, the Shinnecock Inlet is subject to the highest wind speed over the course of the simulation.
One might choose this a priori, as we did in our preliminary numerical experiments. 

We take 25 filtered observations of predicted and observed data using between seven and twelve knots.
Because there are two uncertain parameters $(\lambda_1, \lambda_2)$, we again choose to learn two QoIs.
LUQ chooses a sigmoid kernel, which explains approximately $89\%$ of the variance.
A standard KDE is used to estimate the predicted and observed densities on the learned QoI.
The $\mathbb{E}(\ratio)$ estimate, to two significant digits, is $1.13$.

\begin{figure}
	\centering
	\includegraphics[height=6.5cm]{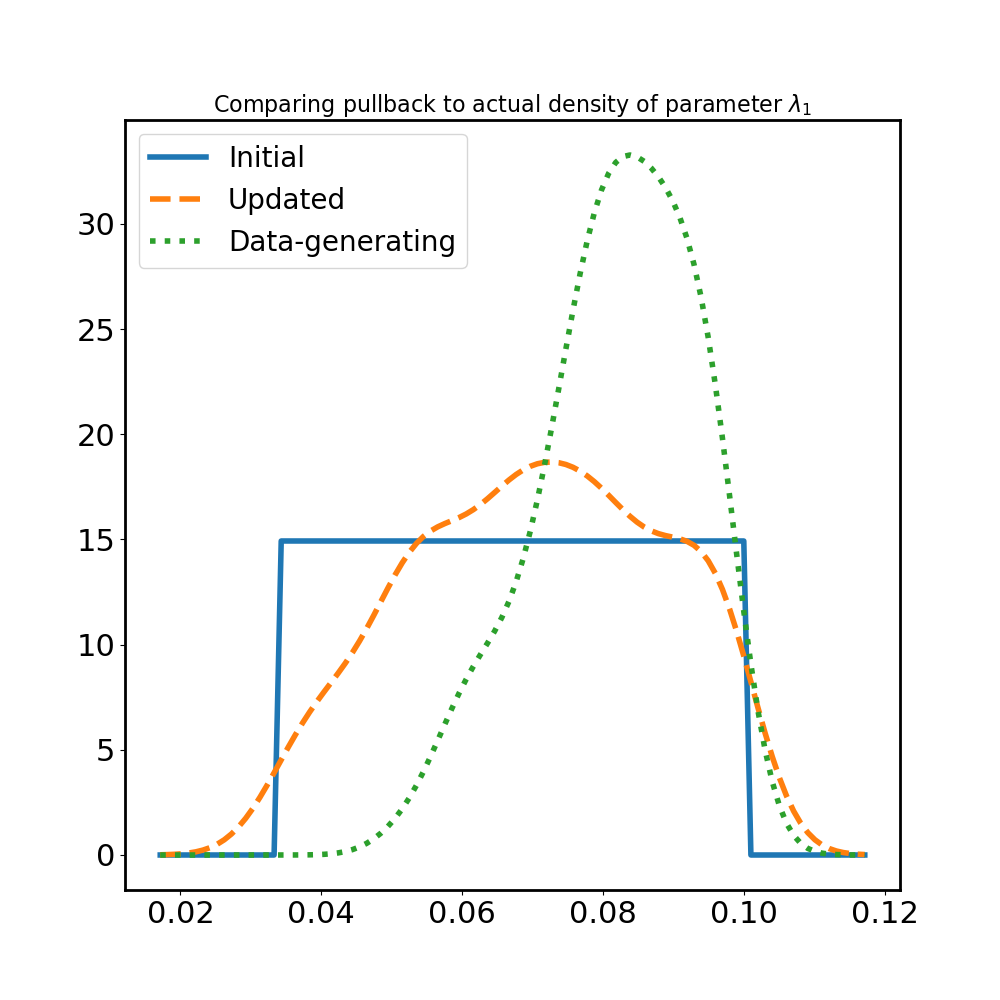}
	\includegraphics[height=6.5cm]{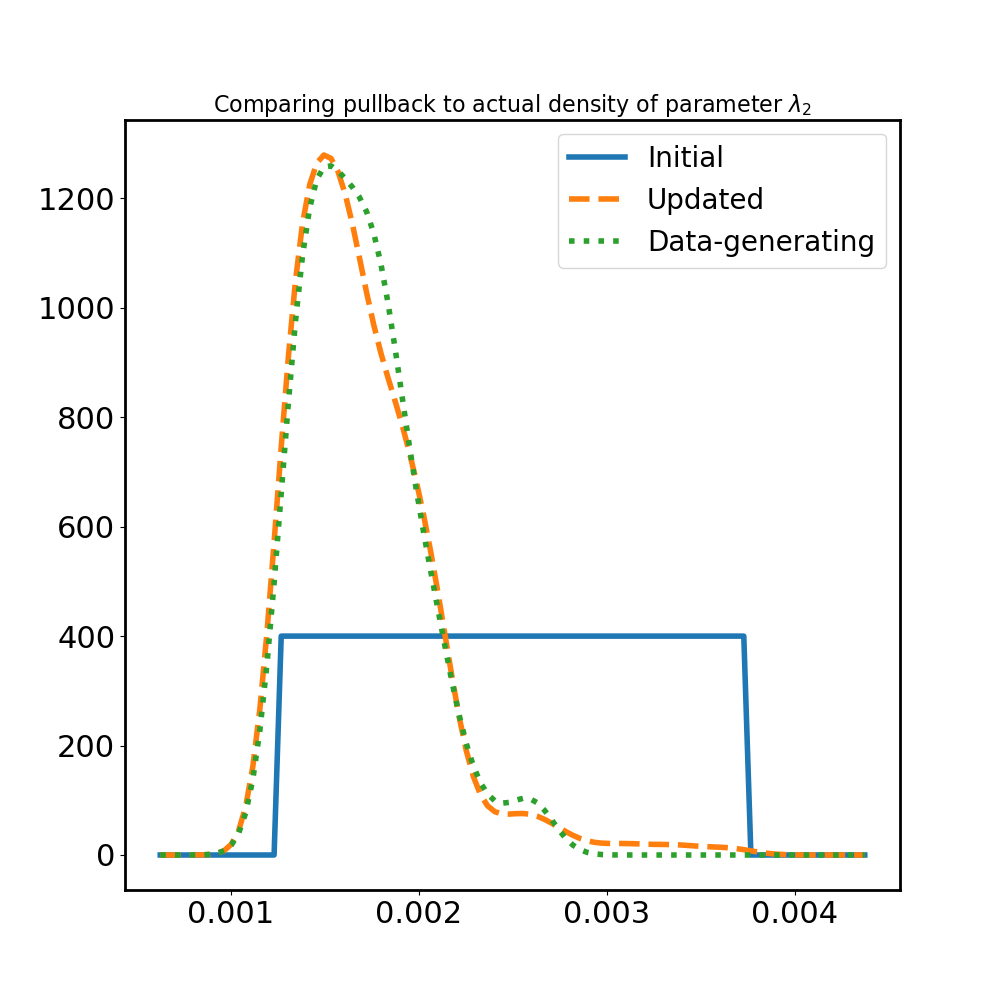}
	\caption{Densities for $\lambda_1$ (left) and $\lambda_2$ (right) for the second numerical experiment presented in Section~\ref{sec:drag}.
	See the caption of Figure~\ref{fig:wd-densities-v1} for interpretation of the plots.}
	\label{fig:wd-densities-v2}
\end{figure}

Figure~\ref{fig:wd-densities-v2} shows several probability densities over $\lambda_1$ (left) and $\lambda_2$ (right).
As before, the solid blue curves are the initial uniform densities, the dashed orange curves show weighted KDEs for the updated densities, while the dotted green curves are standard KDE estimates on the data-generating parameter samples.

The updated density for $\lambda_2$ in Figure~\ref{fig:wd-densities-v2} matches very well with the data-generating density.
In fact, it is an improvement over the previous results where both the mode and tail of the density are better approximated.
This is further quantified in Table~\ref{tab:wd-TV-metrics-2} where we observe a more significant decrease in the total variation distance from the updated marginal density and data-generating marginal density for $\lambda_2$ than in the previous case: dropping from $1.061$ to $0.086$ (second row, Table~\ref{tab:wd-TV-metrics-2}) versus $1.061$ to $0.276$ (second row, Table~\ref{tab:wd-TV-metrics-1}).
However, the updated density for $\lambda_1$ in Figure~\ref{fig:wd-densities-v2} is not significantly different from a KDE estimate of a uniform density and is not a very good estimate of the data-generating density associated with this parameter.


\begin{table}
\centering
\begin{tabular}{|r|r|r|}
\hline
$\| \pi_{\lambda_1}^{init} - \pi_{\lambda_1}^{DG}  \|_{TV}$  &
$\| \pi_{\lambda_1}^{update} - \pi_{\lambda_1}^{DG}  \|_{TV} $ & 
$\| \pi_{\lambda_1}^{DG} - \pi_{\lambda_1}^{DG, exact}  \|_{TV}$ \\
\hline
0.744 & 0.632 & 0.093 \\
\hline
\end{tabular}
\vskip 0.2cm
\begin{tabular}{|r|r|r|}
\hline
$\| \pi_{\lambda_2}^{init} - \pi_{\lambda_2}^{DG}  \|_{TV} $ &
 $ \| \pi_{\lambda_2}^{update} - \pi_{\lambda_2}^{DG}  \|_{TV} $ &
$\| \pi_{\lambda_2}^{DG} - \pi_{\lambda_2}^{DG, exact}  \|_{TV} $ \\
\hline
1.061 & 0.086 & 0.160 \\
\hline
\end{tabular}
\caption{Total variation metrics for the second numerical experiment presented in Section~\ref{sec:drag}, using the time window $T_2$.
See the caption of Table~\ref{tab:wd-TV-metrics-1} for interpretation of the data.}
\label{tab:wd-TV-metrics-2}
\end{table}

In summary, using the new time window $T_2$ results in no significant update in the density for $\lambda_1$ versus excellent agreement for the update in the density for $\lambda_2$.
This, however, is not a failing of the proposed framework for solving the SIP.
Rather, it is a failing in either the choice of experimental design and/or choice of hyperparameters.
The interested reader can readily modify the second half of the provided Python script, publicly available on Archive.org, to use one rather than two QoI; see \ref{app:dependencies} for more detail.
Then, the $\mathbb{E}(\ratio)$ estimate improves dramatically while the plots of density and summary of TV distances (cf., Figure~\ref{fig:wd-densities-v2} and Table~\ref{tab:wd-TV-metrics-2}) are not significantly changed.
One can therefore conclude, with this modification, that the variations in QoI constructed for the SIP over this time window are primarily sensitive to $\lambda_2$. 

A careful physical reasoning provides useful insight into these results as well.
When the wind speeds are higher, the wind drag coefficient $C_d$ is truncated at the value of $\lambda_2$.
This becomes, in essence, the only parameter value ``seen'' by the system over the period of time when the wind speeds are high. 
It is therefore not surprising that data obtained when wind speeds are high exhibit very little sensitivity to the values of $\lambda_1$ and are therefore inadequate in constructing a QoI that updates the initial density on $\lambda_1$ in a meaningful way.

%% file: Conclusions.tex
\section{Conclusions}\label{sec:conclusions}


In this work, we introduce a new conceptual and computational framework, Learning Uncertain Quantities (LUQ), to transform time series data into Quantities of Interest (QoI) that are useful for the observation-consistent solution of stochastic inverse problems arising in the study of dynamical systems.
This provides a rigorous and practical method for uncertainty quantification of dynamical systems using raw streams of time series data.
Key ideas in this framework include data filtering (enabling the analysis of raw, possibly noisy time series), unsupervised learning (to learn dynamics from predicted model output and to classify observations), and feature extraction to determine QoI maps (reducing the dimension of a given time series data set via a clipped kernel PCA).
This extends recent work \cite{BJW18a, BBE11, BES12, BET+14}, in which the QoI map is assumed a priori, to the present setting where instead it must be learned from data.

Numerical results show the efficacy of the proposed framework in solving stochastic inverse problems arising in applications from the life and physical sciences, including the damped harmonic oscillator in Sections~\ref{sec:LUQ}--\ref{sec:inversion}, and a generalization of the Sel'kov model of glycolosis, Burgers' equation, and the depth-averaged shallow water equations (using the ADCIRC coastal ocean model) in Sections~\ref{sec:Hopf}--\ref{sec:drag}, respectively.
In the interest of scientific reproducibility, we also describe our implementation of this new framework, which we have made publicly available on GitHub, as well as provide public access to the data sets and Python scripts for all numerical results presented herein; see \ref{app:dependencies} for details.

Throughout the exposition of the numerical results, we have alluded to an important future direction of research, namely, optimal experimental design.
The LUQ framework permits learning QoI from time series data, but the quality of our results still rely on the sensitivity of the model outputs (time series data) on the model inputs (parameters).
The sensitivity -- or lack thereof -- of the time series data on parameters will in general vary between each parameter, as the parameters change (e.g., due to bifurcations), over different time windows (e.g., due to preasymptotic versus asymptotic regimes in the time series), if time series are extracted from different parts of the domain (e.g., in coastal circulation simulations such as Section~\ref{sec:drag}), etc.
Therefore, an interesting and important direction for future research surrounds how we might ``learn'' an optimal experimental design for a given dynamical system.

Additionally, in the numerical results presented in this work we have used only one time window (TW) per experiment, i.e., one TW for each solution of a stochastic inverse problem.
In contrast, it would be interesting to ``continually learn'' from a time series, for example, by using a sequence of disjoint TWs to solve a sequence of stochastic inverse problems.
We plan to pursue further research in this direction, e.g., establishing under what conditions we might expect to ``converge'' in some sense to a ``steady-state distribution'' on the parameters.

%% file: Software-dependencies.tex
\section{Obtaining software, data, and scripts}\label{app:dependencies}
\verb|LUQ| utilizes several publicly available Python packages that are commonly used for scientific computing (\verb|NumPy| \cite{oliphant2006guide} and \verb|SciPy| \cite{2020SciPy-NMeth}) and machine learning (\verb|scikit-learn| \cite{scikit-learn}).
We suggest using a newer version of Python 3 (Python 3.6 or newer).
Version 1.1 of \verb|LUQ| \cite{mattis_luq} was used in this work.
This version and its required dependencies can be installed using the Python Package Installer (\verb|pip|) by
\begin{verbatim}
pip install git+https://github.com/CU-Denver-UQ/LUQ@v1.1
\end{verbatim}
The most up-to-date version of \verb|LUQ| can be installed by
\begin{verbatim}
pip install git+https://github.com/CU-Denver-UQ/LUQ
\end{verbatim}

This repository also contains contains the scripts that reproduce the figures and table data for the numerical results presented in Sections~\ref{sec:LUQ}--\ref{sec:inversion},~\ref{sec:Hopf}, and~\ref{sec:Shock}.
For the harmonic oscillator problem in Sections~\ref{sec:LUQ}--\ref{sec:inversion} (Figures 1--5, Table 1), see \url{https://github.com/CU-Denver-UQ/LUQ/blob/v1.1/examples/harmonic-oscillator/harmonic_oscillator.py}.
For the Hopf bifurcation problem in Section~\ref{sec:Hopf} (Figures 6--7, Table 2), see \url{https://github.com/CU-Denver-UQ/LUQ/blob/v1.1/examples/selkov/selkov.py}.
Lastly, for the shock problem in Section~\ref{sec:Shock} (Figures 8--11, Table 3), see \url{https://github.com/CU-Denver-UQ/LUQ/blob/v1.1/examples/shock/burgers_shock.py}.

Due to GitHub's file size limitations, the data and scripts for the Shinnecock Inlet problem, presented in Section~\ref{sec:drag}, can be found online at Archive.org (\url{https://archive.org/details/troy-butler-shinnecock-inlet-initial-data}); see \cite{Butler:2020:SID} for further information.
Available at the above URL are three data files, one Python file, and several metadata files; the latter are extraneous for the purposes of this discussion.
All three data files are in the MATLAB ``\verb|.mat|'' file format, which is readily (and automatically) accessible in the provided Python script using a subroutine provided by the aforementioned \verb|SciPy| software package.
The Python script, ``\texttt{Troy Butler - Shinnecock\_Inlet.py}'', utilizes the three data files, together with the \verb|LUQ| package and its dependencies, to reproduce Figures 12--14 and Tables 4--5.